\documentstyle[amstex,amssymb,12pt,leqno]{article}
\title{Local solutions to a free boundary problem\\ for the Willmore functional} 
\author{Roberta Alessandroni and Ernst Kuwert} 
\date{\today}
\makeatletter

\numberwithin{equation}{section}
\newtheorem{theorem}{Theorem}
\newtheorem{corollary}{Corollary}
\newtheorem{lemma}{Lemma}

\newtheorem{proposition}{Proposition}

\def\kasten{\hfill\null\nobreak\hfill \hbox{\vrule\vbox{\hrule width
6pt\vskip6pt\hrule}\vrule} \par\smallskip}

\newcommand{\ve}{\varepsilon}

\newcommand{\R}{{\mathbb R}}

\newcommand{\W}{{\mathcal W}}
\newcommand{\Wbar}{{\overline{{\cal W}}}}

\renewcommand{\S}{{\mathbb S}}
\newcommand{\N}{{\mathbb N}}

\newcommand{\mint}{-\hspace{-1.1em}\int}

\begin{document}
\maketitle
\begin{abstract} \noindent{We consider a free boundary problem for 
the Willmore functional ${\cal W}(f) = \frac{1}{4} \int_\Sigma H^2\,d\mu_f$.
Given a smooth bounded domain $\Omega \subset \R^3$, we construct 
Willmore disks which are critical in the class of surfaces 
meeting $\partial \Omega$ at a right angle along their boundary 
and having small prescribed area. Using rescaling and the implicit 
function theorem, we first obtain constrained solutions with 
prescribed barycenter on $\partial \Omega$. We then study the variation of that 
barycenter.} 
\end{abstract}

\section*{Introduction}
The Willmore energy of an immersed surface $f:\Sigma \to \R^3$ is given by
$$
{\cal W}(f) = \frac{1}{4} \int_\Sigma H^2\,d\mu_f,
$$
for instance ${\cal W}(\S^2) = 4\pi$. Introducing the tracefree 
second fundamental form by decomposing $h = h^\circ + \frac{1}{2} H g$, 
we can write the (scalar) Euler-Lagrange operator as
$$
W[f] = \Delta_g H + |h^\circ|^2 H. 
$$
We study a variational problem for the Willmore energy involving 
a free boundary condition. Let $D = \{z \in \R^2: |z| < 1\}$ and
$\Omega \subset \R^3$ be a given smooth, bounded domain. Putting 
$S = \partial \Omega$ we introduce the class ${\cal M}(S)$ of smooth 
immersions $f:\overline{D} \to \R^3$ meeting $S$ orthogonally from inside 
along $\partial D$, that is 
$$
{\cal M}(S) = \big\{f \in C^\infty(\overline{D},\R^3) \mbox{ immersed}:
f(\partial D) \subset S,\,\frac{\partial f}{\partial \eta} = N^S \circ f
\mbox{ on } \partial D\big\}.
$$
Here $\eta, N^S$ are the interior unit normals of $(D,g)$ and 
$\Omega \subset \R^3$ along the respective boundaries. In the 
(unbounded) special case $\Omega = \R^3_{+}$, the round half-spheres 
$$
\S^2_{+}(a,\lambda)  = a + \lambda \S^2_{+} \quad (a \in \R^2,\,\lambda > 0)
$$
minimize the Willmore energy in the class ${\cal M}(\R^2)$. This 
follows from Simon's monotonicity formula, see \cite{KS04}, after 
reflecting across $\R^2$. In particular, the sphere $\S^2_{+}(a,\lambda)$
minimizes in the smaller class of surfaces $f \in {\cal M}(S)$ 
having the same area ${\cal A}(f) = 2\pi \lambda^2$.
For this variational problem we construct critical points in a general 
domain $\Omega$, provided that the prescribed area is sufficiently 
small. 

{\bf Theorem }{\em\,Let $\Omega \subset \R^3$ be 
a smooth bounded domain, and $S = \partial \Omega$. For each sufficiently small 
$\lambda > 0$ there exist at least two disk-type surfaces $f:D \to \R^3$ 
which are critical points for the Willmore functional restricted to the class
\begin{equation}
\label{eqmlambdaclass}
{\cal M}_\lambda (S) = \{f \in {\cal M}(S): {\cal A}(f) = 2\pi \lambda^2\}.
\end{equation}
Each critical point in ${\cal M}_\lambda(S)$ satisfies, for an appropriate 
$\alpha \in \R$, 
\begin{eqnarray} 
\label{eqcritical}
\Delta_g H + |h^\circ|^2 H & = & \alpha H \hspace{10mm} \mbox{ in }D,\\ 
\label{eqnaturalboundary}
\frac{\partial H}{\partial \eta} + h^S(\nu,\nu) H & = & 0 \hspace{14mm} 
\mbox{ on }\partial D.
\end{eqnarray}}
The proof is based on the implicit function theorem and yields
surfaces which are small, almost-round half-spheres, see 
Corollary \ref{corexistence}. 
We show in addition that as $\lambda \searrow 0$ the constructed surfaces 
concentrate at critical points $a \in S$ of the function $H^S:S \to \R$
(Corollary \ref{corconcentrationpoints}). Reversely, if $a \in S$ is 
a nondegenerate critical point of $H^S$, then there is a local 
family $f_\lambda$ of critical points in ${\cal M}_\lambda(S)$ 
which depends smoothly on $\lambda$ and concentrates at $a$ as 
$\lambda \searrow 0$;  see Theorem \ref{thmlocal} for details.\\
\\
In \cite{Nit93} Nitsche discusses possible boundary conditions for Willmore 
surfaces on grounds of the boundary terms in the first variation 
formula. Palmer proves symmetry and uniqueness for Willmore surfaces
with boundary moving freely on a plane or round sphere \cite{Pal93},
see also Dall'Acqua \cite{Dal12} for related work. It appears that the present 
variational problem involving the class ${\cal M}(S)$ was however not 
considered in the literature. Our main motivation is the conformal 
invariance of the class ${\cal M}(S)$, which should lead to interesting 
compactness and regularity issues. We have verified a reflection principle 
for Willmore surfaces with our boundary condition in the case $\Omega = \R^3_{+}$. 
By the work of Bryant \cite{Bry84}, all disk-type solutions are then
obtained from minimal surfaces with reflectional symmetry, having 
the type of $\S^2$ with finitely many flat ends. Of course one may
also consider the variational problem with other prescribed angles. 
For the one-dimensional Bernoulli
elastic energy and for the Willmore energy under rotational symmetry,
solutions with Dirichlet or Navier type boundary conditions are 
constructed by Deckelnick, Grunau et al., see for instance \cite{DG07,DDG08}. 
Existence and regularity results for Willmore minimizers with prescribed curve
and tangent plane along the boundary were proved by Sch\"atzle \cite{Sch10}. 
A recent paper by  Alexakis and Mazzeo considers properly immersed
surfaces in hyperbolic $3$-space which are (locally) critical points of the 
$L^2$ energy of the second fundamental form. They show that finite energy 
surfaces meet the sphere at infinity at a right angle \cite[Lemma 2.1]{AM12}.\\
\\
To prove the existence result we study the problem on $\R^3_{+}$ 
with respect to pertubations $\tilde{g}$ of the Euclidean background metric. 
On the space of variations of $\S^2_{+}$ respecting the boundary condition, 
the linearized operator has a three-dimensional kernel due to dilations
and translations. We arrive at a solvable problem by prescribing the area 
${\cal A}(f,\tilde{g}) = 2\pi$ and a two-dimensional barycenter $C(f,\tilde{g}) = 0 \in \R^2$.\\
\\
Pulling back the Euclidean metric with a chart near $a \in \partial \Omega$ 
and rescaling yields a perturbed metric $\tilde{g}^{a,\lambda}$ on $\R^3_{+}$.
Solving the constrained problem for $\tilde{g}^{a,\lambda}$ and transforming back, 
we get a three-dimensional family $\phi^{a,\lambda}$ of critical 
points subject to constraints ${\cal A}(\phi^{a,\lambda}) = 2\pi \lambda^2$ 
and $C(\phi^{a,\lambda},S) = a$. In Proposition \ref{propexpansion} we prove 
the expansion
$$
|{\cal W}(\phi^{a,\lambda}) - 2\pi + \pi H^S(a) \lambda| 
\leq C\lambda^2 \quad \mbox{ where }C = C(\Omega). 
$$
In particular $\inf_{f \in {\cal M}(S)} {\cal W}(f) < 2\pi$. This 
indicates that minimizers of ${\cal W}(f)$ without area constraint
are not in the realm of a local approach. In Theorem 
\ref{thmreduction} we show instead the following: for $\lambda \in (0,\lambda_0]$
a constrained solution $\phi^{a,\lambda}$ is critical under prescribed area
${\cal A}(\phi^{a,\lambda}) = 2\pi \lambda^2$ if and only if the
point $a \in S$ is a critical point of the reduced energy function
$$
\bar{{\cal W}}(\cdot,\lambda):S \to \R,\,
\bar{{\cal W}}(a,\lambda) = {\cal W}(\phi^{a,\lambda}).
$$
In consequence we get at least two critical points in ${\cal M}_\lambda(S)$ 
for $\lambda \in (0,\lambda_0(\Omega)]$, as stated in the theorem.\\
\\
In \cite{LamM10,LMS11} Lamm, Metzger and Schulze study a related 
pertubation problem for small spheres in Riemannian manifolds. 
Their solutions are also critical with respect to prescribed area 
and are called {\em of Willmore type}. Another pertubation
result, also in a Riemannian manifold, is by Mondino \cite{Mon10}.\\
\\
There is a corresponding analysis for constant mean curvature
surfaces. The pioneering work is by Ye \cite{Ye91}.
Our approach is close to the work of Pacard and Xu \cite{PX09} and 
also Fall \cite{Fal10,Fal12}. The following difference should 
however be noted: in the CMC case the orthogonality along the
boundary appears as natural boundary condition, whereas here 
it is imposed as a constraint. Our natural boundary condition 
is equation (\ref{eqnaturalboundary}).\\
\\
We now outline the contents of this paper. In Section 1
we compute the space of admissible variations, that is the tangent 
space of ${\cal M}(S)$, and derive the resulting boundary
conditions. One can show that the space ${\cal M}(S)$ is 
a manifold; for the purposes of this paper a graph representation 
of ${\cal M}(S)$ near $\S^2_{+}$ is sufficient (Lemma \ref{lemmamanifold}).
In Section 2 we solve the constrained pertubation problem 
with respect to an arbitrary background Riemannian metric
close to the standard metric. Technically we use a two-step 
procedure where the orthogonality constraint is satisfied 
first, leading to a certain submanifold on which the other 
equations are then solved in the second step, see Lemma 
\ref{lemmaimplicit}.\\
\\
This is applied in Section 3 to the local situation around 
$a \in S$, pulling back and rescaling as indicated above. 
Graph coordinates turn out to be sufficient for this purpose.
We then prove the main results: the expansion of the energy
(Proposition \ref{propexpansion}), the characterization 
of critical points using the reduced energy function (Theorem
\ref{thmreduction}) and finally the existence results 
(Corollary \ref{corexistence} and Theorem \ref{thmlocal}).
In the appendix we review the construction of the two-dimensional
barycenter. 

\section{Constraints and conditions on the boundary}
We start by collecting without proof some variational formulae. 
Let $f:\Sigma \to (M^3,\tilde{g})$ be a compact, smoothly 
immersed surface with boundary $\partial \Sigma$. We denote by $\tilde{D}$ the Levi-Civita connection
on $M$ and by $g=f^*\tilde{g}$ the induced metric on $\Sigma$. We assume that we have a unit normal
$\nu:\Sigma \to TM$ along $f$.

\begin{lemma} \label{lemmavariation} Let $f:\Sigma \times I \to (M^3,\tilde{g})$
be a smooth variation, $0 \in I$, with $\partial_t f = \varphi \nu$ at $t=0$.
Then at $t = 0$ we have the following equations:
\begin{eqnarray}
\label{varderivative}
\tilde{D}_t \partial_k f & = & (\partial_k \varphi) \nu  - \varphi\, g^{ij} h_{jk}\, \partial_i f\\
\label{varmetric}
\partial_t g_{ij} & = & - 2 h_{ij}\, \varphi\\
\label{vararea}
\partial_t (d\mu_g) & = & - H \varphi\, d\mu_g\\
\label{varnormal}
\tilde{D}_t \nu & = & - Df \cdot {\rm grad \varphi}\\
\label{varcurvature}
\partial_t h_{kl} & = & \nabla^2_{kl} \varphi - g^{ij} h_{ik} h_{jl}\, \varphi
+ \tilde{R}(\nu,\partial_k f,\partial_l f,\nu)\, \varphi,\\
\label{varmeancurvature}
\partial_t H & = & \Delta_g \varphi + \big(|h|^2 + {\rm \tilde{Ric}}(\nu,\nu)\big)\, \varphi\\
\label{varchristoffel}
\partial_t \Gamma_{ij}^k & = & 
- g^{kl} \big(\nabla_i h_{jl} + \nabla_j h_{il} - \nabla_l h_{ij}\big) \varphi\\
& & - g^{kl} \big((\partial_i \varphi) h_{jl} + (\partial_j \varphi) h_{il} 
- (\partial_l \varphi) h_{ij} \big).
\end{eqnarray}
In a space of constant curvature $\varkappa$, the curvature terms simplify to
$$
\tilde{R}(\nu,\partial_k f,\partial_l f,\nu) = \varkappa\, g_{kl} 
\quad \mbox{ and } \quad 
{\rm \tilde{Ric}}(\nu,\nu) = 2 \varkappa. 
$$
\end{lemma}

Next we derive the wellknown first variation formula for the Willmore energy.
A version including boundary terms was stated e.g. in \cite{Nit93}.

\begin{theorem} \label{thmwillmorevariation}
For $f:\Sigma \to (M^3,\tilde{g})$, the first variation
of the Willmore energy
in direction of the vector field $\phi = \varphi \nu + Df \cdot \xi$ is
$$ 
\frac{d}{dt} {\cal W}(f)|_{t=0}  = 
\frac{1}{2} \int_\Sigma W(f) \varphi \,d\mu_g 
+ \frac{1}{2} \int_{\partial \Sigma} \omega(\eta)\,ds_g=: \delta{\cal W}(f)\phi,
$$
where $\eta$ is the interior unit normal with respect to $g$, and 
\begin{eqnarray*}
W(f) & = & \Delta H + \big(|h^\circ|^2 + {\rm \tilde{Ric}}(\nu,\nu)\big) H,\\
\omega(\eta) & = & \varphi \frac{\partial H}{\partial \eta} - \frac{\partial \varphi}{\partial \eta} H
- \frac{1}{2} H^2 g(\xi,\eta).
\end{eqnarray*}
\end{theorem}

{\em Proof. } We compute for normal and tangential $\phi$ separately, starting with
the first. In normal coordinates for $t = 0$ we get from Lemma \ref{lemmavariation}
\begin{eqnarray*}
\frac{d}{dt} {\cal W}(f) & = & 
\frac{1}{2} \int_\Sigma \frac{\partial H}{\partial t} H \,d\mu_g
+ \frac{1}{4} \int_\Sigma H^2\,\frac{\partial}{\partial t}\,d\mu_g\\
& = & \frac{1}{2} \int_\Sigma \Big(
\Delta \varphi + \big(|h|^2 + {\rm \tilde{Ric}}(\nu,\nu)\big) \varphi \Big) H \,d\mu_g
- \frac{1}{4} \int_\Sigma H^3\varphi \,d\mu_g.
\end{eqnarray*}
Using $|h|^2 = |h^\circ|^2 + \frac{1}{2} H^2$ and integrating by parts yields
\begin{eqnarray*}
\frac{d}{dt} {\cal W}(f) & = & \frac{1}{2} 
\int_\Sigma \big(\Delta H + (|h^\circ|^2 + {\rm \tilde{Ric}}(\nu,\nu)) H\big) \varphi\,d\mu_g\\
& & + \frac{1}{2} \int_{\partial \Sigma} 
\Big(\varphi \frac{\partial H}{\partial \eta} - \frac{\partial \varphi}{\partial \eta} H\Big)\,ds_g. 
\end{eqnarray*}
This proves the claim in the case when $\phi$ is normal. Now consider a variation
of the form $f \circ \varphi_t$, where $\varphi_t$ is the flow of the vector
field $\xi$. For $Q\subset \!\! \subset \Sigma$ we get by invariance
with respect to reparametrizations
\begin{eqnarray*}
{\cal W}(f \circ \varphi_t,Q) & = & {\cal W}(f,\varphi_t(Q))\\
& = & \frac{1}{4} \int_{Q_t} H(y)^2\,d\mu_g(y)\\
& = & \frac{1}{4} \int_{Q} H(\varphi_t(x))^2\,J\varphi_t(x)\,d\mu_g(x),
\end{eqnarray*}
where $J\varphi_t(x)$ is the Jacobian. Differentiating at $t = 0$ we get
\begin{eqnarray*}
\frac{d}{dt} {\cal W}(f \circ \varphi_t,Q) & = & \frac{1}{4}  \int_{Q}
\big(\partial_{\xi} H^2 + H^2 {\rm div}_g \xi\big)\,d\mu_g\\
& = & \frac{1}{4}  \int_{Q} {\rm div}_g (H^2 \xi)\,d\mu_g\\
& = & - \frac{1}{4}  \int_{\partial Q} H^2 g(\xi,\eta) \,ds_g.
\end{eqnarray*}
Since $\omega(\eta) = - \frac{1}{2} H^2 g(\xi,\eta)$ for $\phi = Df \cdot \xi$ (i.e.
$\varphi \equiv 0$), the formula is proved for all $\phi$. \kasten

Now let $\Omega \subset \R^3$ be a domain with smooth boundary. We put 
$S = \partial \Omega$ and denote by $N^S:S \to \S^2$ the interior unit 
normal. Then for a smooth compact surface $\overline{\Sigma} = \Sigma \cup \partial \Sigma$
we consider the class of immersions 
\begin{equation}
\label{eqclassM}
{\cal M}(S) = \big\{f \in C^\infty(\overline{\Sigma},\R^3) \mbox{ immersed: }
f(\partial \Sigma) \subset S,\,\frac{\partial f}{\partial \eta} = N^S \circ f\big\}. 
\end{equation}
Let $h$ and $h^S$ be the second fundamental forms of $f$ and $S$, respectively. 
We calculate, using that $\tilde{D}^2 f(\tau,\tau) = h(\tau,\tau) \nu$ is normal 
to $N^S \circ f$, for $\tau$ the unit tangent along $\partial \Sigma$, 
$$
0 = \frac{\partial}{\partial \tau} \tilde{g}\left(\frac{\partial f}{\partial \tau}, N^S \circ f\right)
= \tilde{g}\left( Df \cdot \nabla_\tau \tau, N^S \circ f\right)
+ \tilde{g}\left( \frac{\partial f}{\partial \tau}, (DN^S) \circ f \cdot \frac{\partial f}{\partial \tau} \right).
$$
The geodesic curvature of $\partial \Sigma$ with respect to 
the induced metric $g$ is defined by
$$
\nabla_\tau \tau = \varkappa_g \eta \quad \Leftrightarrow \quad 
\nabla_\tau \eta = - \varkappa_g \tau.
$$
Thus $\varkappa_g = +1$ for the standard disk. We have 
\begin{equation}
\label{eqconstraint4}
\varkappa_g = h^S\big(\frac{\partial f}{\partial \tau},\frac{\partial f}{\partial \tau}\big).
\end{equation}
Taking the derivative of $\tilde{g}(\nu,N^S) = 0$ in the direction of $\tau$ yields
\begin{equation}
h(\tau,\eta) + h^S(\nu,\frac{\partial f}{\partial \tau}) = 0. 
\label{eqconstraint3}
\end{equation}
A further tangential derivative implies
\begin{equation}
\nabla_\tau h(\tau,\eta) + \varkappa_g (h(\eta,\eta)-h(\tau,\tau))
+ \frac{\partial}{\partial \tau} \big[h^S(\nu,\frac{\partial f}{\partial \tau})\big] = 0.\label{eqconstraint5}
\end{equation}
Next we linearize the constraints. Let $f = f(p,t) \in {\cal M}(S)$ and put
$$
\frac{\partial f}{\partial t}|_{t=0} =\phi = \varphi \nu + Df \cdot \xi.
$$
Differentiating the equation $f(\partial \Sigma,t) \in S$ yields
\begin{equation}
\label{eqlinboundary1}
0 = \tilde{g}\left( \phi, N^S \circ f\right) = \tilde{g}\left( Df \cdot \xi,Df \cdot \eta \right)
= g(\xi,\eta) \quad \mbox{ along } \partial \Sigma.
\end{equation}
For the variation of the normal we have the standard formula
\begin{equation}
\label{eqvariationnormal}
\frac{\tilde{D} \nu}{\partial t} = Df \cdot \big(-{\rm grad}_g \varphi + W \xi\big) 
\quad \mbox{ on } \Sigma, 
\end{equation}
where $W$ is the Weingarten map given by $h(X,Y) = - g(WX,Y)$
or $\tilde{D}\nu = Df \cdot W$. The first variation of the orthogonality 
relation gives
$$
0 = \frac{\partial}{\partial t} \tilde{g}\left( \nu,N^S \circ f\right)
= \tilde{g}\left( Df \cdot \big(-{\rm grad}_g \varphi + W \xi), N^S \circ f\right)
+ \tilde{g}\left( \nu, (W^S \circ f) \phi \right).
$$
In this calculation we used $f(\partial \Sigma,t) \subset S$ so that
$\phi \in T_f S$  and $(W^S \circ f) \phi$ makes sense. 
Now from $N^S \circ f = \frac{\partial f}{\partial \eta}$
we have
$$
\tilde{g}\left( Df \cdot \big(-{\rm grad}_g \varphi + W \xi), N^S \circ f\right)
= g(-{\rm grad}_g \varphi + W \xi,\eta)
= - \frac{\partial \varphi}{\partial \eta} - h(\xi,\eta).
$$
Using further $\tilde{g}\left( \nu,(W^S \circ f)\phi \right) = - h^S(\nu,\varphi \nu + Df \cdot \xi)$
we arrive at the following two linearized equations, for the variation vectorfield
$\phi = \varphi \nu + Df \cdot \xi$,
\begin{eqnarray}
\label{eqlinboundary2}
g(\xi,\eta) & = & 0 \quad \mbox{ on } \partial \Sigma,\\
\label{eqlinboundary3}
\frac{\partial \varphi}{\partial \eta} + h(\xi,\eta) 
+ \varphi h^S(\nu,\nu) + h^S(\nu,Df \cdot \xi) & = & 0 \quad \mbox{ on }\partial \Sigma.
\end{eqnarray}
Equation (\ref{eqlinboundary2}) holds if and only if $\xi = \mu \tau$ for some
function $\mu$ on $\partial \Sigma$. Then (\ref{eqlinboundary3}) simplifies using
(\ref{eqconstraint3})  and we are left with
\begin{equation}
\label{eqlinboundary4}
\frac{\partial \varphi}{\partial \eta} + \varphi h^S(\nu,\nu) = 0
\quad \mbox{ where } \phi = \varphi \nu + \mu \frac{\partial f}{\partial \tau} 
\mbox{ on } \partial \Sigma. 
\end{equation}
The variation vector fields $\phi$ with (\ref{eqlinboundary4}) are called 
{\em admissible for $f$} and are denoted by $T_f{\cal M}(S)$. Any function $\varphi$ given on 
$\partial \Sigma$ admits an extension to $\Sigma$ such that (\ref{eqlinboundary4}) holds,
and for any $\mu$ on $\partial \Sigma$ there exists a vector field $\xi$ on $\Sigma$
such that $\xi|_{\partial \Sigma} = \mu \tau$. Then the variation
$\phi = \varphi \nu + Df \cdot \xi$ is admissible.\\
\\
Now assume that $f \in {\cal M}(S)$ satisfies 
\begin{equation}
\label{eqcriticalpoint}
\delta{\cal W}(f) \phi = 0 \quad 
\quad \mbox{ for all } \phi \in T_f{\cal M}(S).
\end{equation}
Then clearly $W(f) = 0$, and the definition
of $T_f{\cal M}$ as in (\ref{eqlinboundary4}) implies further 
$$
0 = \frac{1}{2} \int_{\partial \Sigma}
\varphi \Big(\frac{\partial H}{\partial \eta} + H h^S(\nu,\nu)\Big)\,ds_g
\quad \mbox{ for all } \varphi \in C^\infty(\partial \Sigma).
$$
So we arrive at the two boundary conditions
\begin{eqnarray}
\label{eqboundary1}
\tilde{g}\left( \nu, N^S \circ f \right) & = & 0 \quad \mbox{ on } \partial \Sigma,\\
\label{eqboundary2}
\frac{\partial H}{\partial \eta} + H h^S(\nu,\nu) & = & 0 \quad \mbox{ on } \partial \Sigma.
\end{eqnarray}
This paper studies a perturbed boundary value problem with respect 
to Riemannian metrics $\tilde{g}$ which are close to the Euclidean 
metric $\delta$, aiming at immersions close to the standard
$\S^2_{+}$. We now collect some formulae for radial graphs
$$
f:{\mathbb S}^2_{+} \to \R^3,\,f(\omega) = (1+w(\omega)) \omega.
$$
For a tangent vector $\tau \in T_\omega{\mathbb S}^2$ we have
$$
\partial_\tau f(\omega) = (1+w(\omega))\tau + (\partial_{\tau} w)(\omega) \omega.
$$
In an orthonormal frame $\tau_1,\tau_2$ on ${\mathbb S}^2$ the metric 
$g = f^\ast \tilde{g}$ is given by 
\begin{eqnarray*}
g(\tau_\alpha,\tau_\beta) & = & (1+w)^2 \tilde{g}(\tau_\alpha,\tau_\beta)\\
&& + (1+w)(\partial_{\tau_\alpha} w) \tilde{g}(\omega,\tau_\beta) 
+ (1+w)(\partial_{\tau_\beta} w) \tilde{g}(\omega,\tau_\alpha)\\ 
&& + (\partial_{\tau_\alpha} w) (\partial_{\tau_\beta} w) \tilde{g}(\omega,\omega).
\end{eqnarray*}
Here $\tilde{g}$ is always evaluated at $f(\omega)$. The area of $f$ with respect to 
$\tilde{g}$ is 
$$
{\cal A}(f,\tilde{g}) = \int_{{\mathbb S}^2_{+}} \sqrt{\det g(\tau_\alpha,\tau_\beta)}\,d\mu_{{\mathbb S}^2}.
$$
Let $\tilde{\nu}_{\R^2}$ be the upper unit normal along $\R^2$ with respect to $\tilde{g}$. 
We compute 
$$
{\rm grad}_{\tilde{g}}\,x^3 = \sum_{i,j =1}^3 \tilde{g}^{ij} \partial_i x^3 e_j
= \sum_{j=1}^3 \tilde{g}^{3j} e_j.
$$
Further $\tilde{g}({\rm grad}_{\tilde{g}}\,x^3,{\rm grad}_{\tilde{g}}\,x^3) = 
\tilde{g}_{jk} \tilde{g}^{3j} \tilde{g}^{3k} = \tilde{g}^{33}$. Thus we have 
$$
\tilde{\nu}_{\R^2} = \frac{1}{\sqrt{\tilde{g}^{33}}} \sum_{j=1}^3 \tilde{g}^{3j} e_j.
$$
Now let $\nu:{\mathbb S}^2_{+} \to \R^3$, $\nu = \nu[f,\tilde{g}]$, be the unit normal along 
$f$ with respect to $\tilde{g}$, such that $\nu(\omega) = - \omega$ for $u = 0$, 
$\tilde{g} = \delta$. Then 
$$
\tilde{g}(\nu,\tilde{\nu}_{\R^2}) =  
\frac{1}{\sqrt{\tilde{g}^{33}}} \tilde{g}(\nu, {\rm grad}_{\tilde{g}}\,x^3))
= \frac{1}{\sqrt{\tilde{g}^{33}}} \langle \nu,e_3 \rangle. 
$$
With respect to $\tilde{g}(f(\omega))$, the component of $\omega$ which is tangential along $f$ is 
$$
\omega^\top = g^{\alpha \beta} \tilde{g}(\omega,\partial_{\tau_\alpha} f) \partial_{\tau_\beta} f.
$$
Here $\tilde{g}$ is always evaluated at $f(\omega)$. Then $\omega^\perp = \omega-\omega^\top$ 
has the norm 
$$
\tilde{g}(\omega^\perp,\omega^\perp) 
= \tilde{g}(\omega,\omega-\omega^\top)
= \tilde{g}(\omega,\omega) 
- g^{\alpha \beta} \tilde{g}(\omega,\partial_{\tau_\alpha} f) \tilde{g}(\omega,\partial_{\tau_\beta} f).
$$
Dividing we obtain the formula
$$
\nu(\omega) = - \frac{\omega - g^{\alpha \beta} \tilde{g}(\omega,\partial_{\tau_\alpha} f) \partial_{\tau_\beta} f}
{\sqrt{\tilde{g}(\omega,\omega) 
- g^{\alpha \beta} \tilde{g}(\omega,\partial_{\tau_\alpha} f) \tilde{g}(\omega,\partial_{\tau_\beta} f)}}.
$$
The following two lemmas show that the constraint of orthogonality is nondegenerate at the 
standard $\S^2_{+}$. 

\begin{lemma} We have $W^{2,2}(\S^2_{+}) = X_0 \oplus Y_0$ as topological direct sum, where
\begin{eqnarray*}
X_0 & = & \{u \in W^{2,2}(\S^2_{+}): \frac{\partial u}{\partial \eta} = 0 
\quad \mbox{ on } \partial \S^2_+ \}\\
Y_0 & = & \{v \in W^{2,2}(\S^2_{+}): \Delta_{\S^2} v = {\rm const.}\mbox{ on } \S^2_+,\, \int_{\S^2_{+}} v \,d\mu_g= 0\}.
\end{eqnarray*}
Moreover $C^{k,\alpha}(\S^2_{+}) = (X_0 \cap C^{k,\alpha}(\S^2_{+})) \oplus (Y_0 \cap C^{k,\alpha}(\S^2_{+}))$
for any $k \geq 2$, $\alpha \in (0,1)$.
\end{lemma} 

{\em Proof. }$X_0$ and $Y_0$ are closed subspaces of $W^{2,2}(\S^2_{+})$ with $X_0 \cap Y_0 = \{0\}$. 
Any $w \in W^{2,2}(\S^2_{+})$ decomposes uniquely as $w = u + v$, where $u \in X_0$, $v \in Y_0$ are 
chosen with 
\begin{eqnarray*}
\Delta_{\S^2} v & = & - \frac{1}{2\pi} \int_{\partial \S^2_{+}} 
\frac{\partial w}{\partial \eta}\,ds_g \quad \mbox{ in }\S^2_{+},\quad 
\frac{\partial v}{\partial \eta} = \frac{\partial w}{\partial \eta} \mbox{ on }\partial \S^2_{+}\\
u & =  & w - v. 
\end{eqnarray*}
Using Sobolev trace theory \cite{Ada75,Mor66} we have the a priori estimates
$$
\|u\|_{W^{2,2}(\S^2_{+})} + \|v\|_{W^{2,2}(\S^2_{+})} \leq C\,\|w\|_{W^{2,2}(\S^2_{+})}.
$$
Therefore the map $X_0 \oplus Y_0 \to W^{2,2}(\S^2_{+})$, $(u,v) \mapsto u+v$,
is an isomorphism of Banach spaces. Moreover by Schauder regularity \cite{GT83,Mor66} for the 
Neumann problem 
$$
\|u\|_{C^{k,\alpha}(\S^2_{+})} + \|v\|_{C^{k,\alpha}(\S^2_{+})} \leq C\,\|w\|_{C^{k,\alpha}(\S^2_{+})}.
$$
This proves the second statement. \kasten

In the following calculations we assume the background metric $\tilde{g}$ to be 
given on the cylinder $Z_2 = D_2(0) \times [-2,2]$, which compactly contains 
the ball $B_1(0)$.

\begin{lemma} \label{lemmamanifold} Let $\nu = \nu[w,\tilde{g}]$ denote the unit normal of 
the graph of $w \in C^{k,\alpha}(\S^2_{+})$ with respect to the Riemannian metric 
$\tilde{g} \in C^l(Z_2,\R^{3 \times 3}_{{\rm sym}})$. For $1 \leq k \leq l$ 
the map
$$
B[w,\tilde{g}] = \tilde{g}(\nu,\tilde{\nu}_{\R^2}) 
= \frac{1}{\sqrt{\tilde{g}^{33}}} \langle \nu,e_3 \rangle|_{\partial \S^2_{+}} 
\in C^{k-1,\alpha}(\partial\S^2_+)
$$
is well-defined and of class $C^{l-k}$. 
For $2 \leq k < l$ there exist open neighborhoods $U \subset X_0 \cap C^{k,\alpha}(\S^2_{+})$,
$V \subset Y_0 \cap C^{k,\alpha}(\S^2_{+})$ and 
$G \subset C^l(Z_2,\R^{3 \times 3}_{{\rm sym}})$ of 
$u \equiv 0$, $v \equiv 0$ and $\tilde{g} \equiv \delta$, and a $C^{l-k}$ map
$\Psi:U \times G \to V$ such that for all $u \in U$, $v \in V$, $\tilde{g} \in G$
we have 
$$
B[u+v,\tilde{g}] = 0 \quad \Leftrightarrow \quad v = \Psi[u,\tilde{g}].
$$
We have $D_u\Psi[0,\delta] = 0$, and 
$h = D_{\tilde{g}}\Psi[0,\delta] q \in Y_0 \cap C^{k,\alpha}(\S^2_{+})$ is 
the unique solution of 
\begin{equation}
\label{eqmanifold1}
- \Delta_{\S^2} h = \frac{1}{2\pi} \int_{\partial \S^2_{+}} q(\nu,e_3)\,ds \mbox{ in } \S^2_{+},\quad 
\frac{\partial h}{\partial \eta} = q(\nu,e_3) \mbox{ on } \partial \S^2_{+}.
\end{equation}
\end{lemma}

{\em Proof. }The map $B:C^{k,\alpha}(\S^2_{+}) \times C^l(Z_2,\R^{3 \times 3}_{{\rm sym}}) 
\to C^{k-1,\alpha}(\partial \S^2_{+})$ is well-defined and of class $C^{l-k}$ near $w \equiv 0$, 
$\tilde{g} \equiv \delta$, and has the derivative
$$
D_w B[0,\delta] \varphi = -\frac{\partial \varphi}{\partial \eta},
$$
thus ${\rm ker\,}D_w B[0,\delta] = X_0 \cap C^{k,\alpha}(\S^2_{+})$. The operator
$D_w B[0,\delta]|_{Y_0}:Y_0 \cap C^{k,\alpha}(\S^2_{+}) \to C^{k-1,\alpha}(\partial \S^2_{+})$
is an isomorphism: for any $\beta \in C^{k-1,\alpha}(\partial \S^2_{+})$
there is a unique $v \in Y_0 \cap C^{k,\alpha}(\S^2_{+})$ with
$D_w B[0,\delta] v = \beta$, in other words
$$
- \Delta_{\S^2}  v= \frac{1}{2\pi} \int_{\partial \S^2_{+}} \beta\,ds_g, \quad \int_{\S^2_{+}} v \, d\mu_g= 0,\quad
\frac{\partial v}{\partial \eta} = \beta.
$$
Moreover that solution $v$ satisfies the estimate
$\|v\|_{C^{k,\alpha}(\S^2_{+})} \leq C \|\beta\|_{C^{k-1,\alpha}(\partial \S^2_{+})}$,
which means that $D_w B[0,\delta]|_{Y_0}$ has a bounded inverse. Existence and 
uniqueness of $\Psi[u,\tilde{g}]$ follows from the implicit function theorem. 
Now $\Psi[0,\delta] = 0$, and we have for any $\varphi \in X_0 \cap C^{k,\alpha}(\S^2_{+})$ 
$$
0 = \frac{d}{dt} B\big[t\varphi + \Psi(t\varphi,\delta),\delta\big]|_{t=0} = 
\underbrace{D_w B[0,\delta] \varphi}_{=0}  + D_w B[0,\delta] D_u\Psi[0,\delta] \varphi.
$$
This shows $D_u\Psi[0,\delta] = 0$. We have 
further for $\tilde{g} = \delta + tq$ and $\nu = \nu[0,\tilde{g}]$
$$
D_{\tilde{g}} B[0,\delta] \cdot q
= \langle \frac{\partial \nu}{\partial t}|_{t=0}, e_3 \rangle
= \frac{\partial}{\partial t} \underbrace{\tilde{g}(\nu,e_3)}_{=0}|_{t=0} 
- \frac{\partial \tilde{g}}{\partial t}(\nu,e_3)|_{t=0}
= - q(\nu,e_3),
$$
which yields the remaining claim, namely
$$
0 = \frac{d}{dt} B\big[\Psi[0,\delta + tq], \delta + tq\big]|_{t=0} = 
D_w B[0,\delta] D_{\tilde{g}}\Psi[0,\delta] \cdot q - q(\nu,e_3). 
$$
\kasten

\section{The Riemannian pertubation problem}
\setcounter{equation}{0}
Using reflection and Simon's monotonicity formula, it is easy to see that the 
standard half-sphere $\S^2_{+}$ minimizes the Willmore functional among
surfaces meeting $\R^2$ orthogonally. One might hope to get corresponding 
Willmore surfaces for perturbed background metrics $\tilde{g}$ using the
implicit function theorem. However the linearized problem has a kernel $K_0$. 
For any $\lambda > 0$ the dilated sphere $\lambda \S^2_{+}$ also minimizes,
and is represented as graph of $w^\lambda(\omega) \equiv \lambda -1$ 
over $\S^2_{+}$. Hence $K_0$ contains the function 
$$
\frac{\partial}{\partial \lambda} w^\lambda|_{\lambda = 1} \equiv 1
$$
Likewise for any $a \in \R^2$, $|a| < 1$, the translated
halfspheres $\S^2_{+}(a)$ admit the graph representations
$w^a(\omega) = \langle \omega,a \rangle - 1 + \sqrt{1-|a|^2 + \langle \omega,a \rangle^2}$
over $\S^2_{+}$, hence $K_0$ also contains the functions 
$$
\frac{\partial}{\partial \ve} w^{\ve a}(\omega)|_{\ve = 0} = \langle \omega,a \rangle.
$$
We get a solvable problem by prescribing the Riemannian area and two-dimensional 
barycenter. For these constrained solutions the Willmore operator is in the 
space $K(\tilde{g})$ spanned by the $L^2$ gradients of the constraints,
and we have $K(\delta) = K_0$. In the next section we will study the Willmore 
energy as a function on the manifold of constrained solutions.

\begin{lemma} \label{lemmaiso1} Let $K_0 = {\rm Span\,}\{1,\langle \omega,e_1 \rangle,\langle \omega,e_2 \rangle\}$,
and define the Hilbert space 
$$
W^{2,2}_{0,\perp}(\S^2_{+}) = \{u \in W^{2,2}(\S^2_{+}): \frac{\partial u}{\partial \eta} = 0 
\mbox{ on } \partial \S^2_{+},\, u \perp_{L^2} K_0\}.
$$
Then the linear operator 
$$
L:W^{2,2}_{0,\perp}(\S^2_{+})  \to W^{2,2}_{0,\perp}(\S^2_{+})', \langle Lu,v \rangle = \int_{\S^2_{+}} 
\Big(\Delta_{\S^2} u\, \Delta_{\S^2} v - 2 \langle \nabla u,\nabla v \rangle\Big)\,d\mu_{\S^2} 
$$
is an isomorphism.
\end{lemma}
{\em Proof. }Let $E_k \subset L^2(\S^2)$, $k \in \N_0$, be the space of even eigenfunctions of $-\Delta_{\S^2}$
on the $2$-sphere, with eigenvalue $\lambda_k = k (k+1)$ (even means $u(x,z) = u(x,-z)$). We have 
\begin{eqnarray*}
2 \langle L u_k,u_l\rangle & = &
\int_{\S^2} \big(\Delta_{\S^2} u_k \Delta_{\S^2} u_l- 2 \langle \nabla u_k,\nabla u_l \rangle\big)\,d\mu_{\S^2}\\
& = & \int_{\S^2} \Delta_{\S^2} u_k (\Delta_{\S^2} + 2) u_l\,d\mu_{\S^2}\\
& = & \lambda_k (\lambda_{l}-2) \langle u_k, u_l \rangle_{L^2(\S^2)}. 
\end{eqnarray*}
Now $\lambda_k \geq 6$ for $k \geq 2$, thus for a finite sum $u = \sum_{k=2}^N u_k$ we see
$$
\int_{\S^2} \big(\left(\Delta_{\S^2} u\right)^2 + u^2\big)\,d\mu_{\S^2}
= \sum_{k = 2}^N (\lambda_k^2 +1) \|u_k\|^2_{L^2(\S^2)} 
\leq \frac{37}{24} \sum_{k = 2}^N \lambda_k(\lambda_k-2) \|u_k\|^2_{L^2(\S^2)}
=  \frac{37}{12} \langle Lu, u \rangle.
$$
Applying the Bochner Formula on $\S^2$ we conclude that
$$
\int_{\S^2}\left(\vert\nabla^2 u\vert^2+\vert \nabla u\vert^2+ u^2 \right)\,d\mu_{\S^2}
=\int_{\S^2}\left(\left(\Delta_{\S^2} u\right)^2 +u^2 \right)\,d\mu_{\S^2}\\
\leq \frac{37}{12}\langle Lu,u \rangle.
$$
Extending functions $u \in W^{2,2}_{0,\perp}(\S^2_{+})$ by even reflection across $\partial \S^2_{+}$ yields
$W^{2,2}$-functions on the sphere. It is then easy to see that the algebraic sum $\bigoplus_{k=2}^\infty E_k$ 
is $W^{2,2}$-dense in $W^{2,2}_{0,\perp}(\S^2_{+})$. The coercivity of $L$ and hence the 
claim of the lemma follows. \kasten

\begin{lemma} \label{lemmaiso2} For $k \geq 4$ and $\alpha \in (0,1)$ the linear operator 
$$
{\cal L}:C^{k,\alpha}_{0,\perp}(\S^2_{+}) \to C^{k-4,\alpha}_{\perp}(\S^2_{+}) \times C^{k-3,\alpha}(\partial \S^2_{+}),\quad 
{\cal L} u = \Big(\Delta_{\S^2} (\Delta_{\S^2} + 2) u, \frac{\partial (\Delta_{\S^2} u)}{\partial \eta}\Big),
$$
is an isomorphism. 
\end{lemma}

{\em Proof. }We have a commutative diagram 
$$
\begin{array}{ccc} C^{k,\alpha}_{0,\perp}(\S^2_{+}) 
& \stackrel{{\cal L}}{\longrightarrow} & C^{k-4,\alpha}_{\perp}(\S^2_{+}) \times C^{k-3,\alpha}(\partial \S^2_+)\\
\bigcap & & \bigcap\\
W^{2,2}_{0,\perp}(\S^2_{+}) & \stackrel{L}{\longrightarrow} & W^{2,2}_{0,\perp}(\S^2_{+})'
\end{array}
$$
Here for $(f_0,f_1) \in C^{k-4,\alpha}_{\perp}(\S^2_{+}) \times C^{k-3,\alpha}(\partial \S^2_+)$ 
the inclusion on the right is given by 
$$
\Lambda(v) = \int_{\S^2_{+}} f_0 v\,d\mu_{\S^2} + \int_{\partial \S^2_{+}} f_1 v\,ds.
$$
The injectivity of ${\cal L}$ follows from Lemma \ref{lemmaiso1}. Moreover 
for given $(f_0,f_1) \in C^{k-4,\alpha}_{\perp}(\S^2_{+}) \times C^{k-3,\alpha}(\partial \S^2_+)$ 
there exists $u \in W^{2,2}_{0,\perp}(\S^2_{+})$ such that 
$$
\int_{\S^2_{+}} \big(\Delta_{\S^2} u\, \Delta_{\S^2} v - 2 \langle \nabla u,\nabla v \rangle\big)\,d\mu_{\S^2} 
= \int_{\S^2_{+}} f_0 v\,d\mu_{\S^2} + \int_{\partial \S^2_{+}} f_1 v\,ds \quad 
\mbox{ for all } v \in W^{2,2}_{0,\perp}(\S^2_{+}).
$$
This means that $u \in W^{2,2}_{0,\perp}(\S^2_{+})$ solves the equations 
$$
\Delta_{\S^2} (\Delta_{\S^2} + 2) u = f_0 \mbox{ in } \S^2_{+},\quad
\frac{\partial (\Delta_{\S^2} u)}{\partial \eta} = f_1 \mbox{ on } \partial \S^2_{+}.
$$
In fact, integrating by parts for functions $u,v \in C^4(\overline{\S^2_{+}})$ yields
\begin{eqnarray*} 
&& \int_{\S^2_{+}} \Delta_{\S^2} (\Delta_{\S^2} + 2) u \cdot v\,d\mu_{\S^2}\\ 
& = & \int_{\S^2_{+}} \Big({\rm div\,}[\nabla (\Delta_{\S^2} u + 2u) \cdot v] \,d\mu_{\S^2}
- \int_{\S^2_{+}} \langle \nabla (\Delta_{\S^2} u + 2u),\nabla v\rangle\Big)\,d\mu_{\S^2}\\
& = & \int_{\S^2_{+}} \Big({\rm div\,}[\nabla (\Delta_{\S^2} u + 2u) \cdot v] \,d\mu_{\S^2}
- \int_{\S^2_{+}} {\rm div\,}[(\Delta_{\S^2} u + 2u) \nabla v]\,d\mu_{\S^2}\\
&& + \int_{\S^2_{+}} (\Delta_{\S^2} u + 2u) \Delta_{\S^2} v \,d\mu_{\S^2}\\
& = & \int_{\S^2_{+}} \Big(\Delta_{\S^2} u \Delta_{\S^2} v - 2 \langle \nabla u,\nabla v \rangle\Big) \,d\mu_{\S^2}
- \int_{\partial \S^2_{+}} u \frac{\partial v}{\partial \eta}\,ds\\
&& - \int_{\partial \S^2_{+}} \frac{\partial (\Delta_{\S^2} u + 2u)}{\partial \eta} \cdot v\,ds
+ \int_{\partial \S^2_{+}} (\Delta_{\S^2}u+2u) \frac{\partial v}{\partial \eta}\,ds.
\end{eqnarray*}
Schauder theory, see \cite{ADN59,Mor66}, implies $u \in C^{k,\alpha}_{0,\perp}(\S^2_{+})$ and 
$$
\|u\|_{C^{k,\alpha}(\S^2_{+})} \leq C 
\big(\|f_0\|_{C^{k-4,\alpha}(\S^2_{+})} + \|f_1\|_{C^{k-3,\alpha}(\partial \S^2_{+})}\big).
$$
This proves the lemma. \kasten

Now consider on $Z_2 = \{(x,z) \in \R^2 \times \R: |x|, |z| < 2\}$ a given 
Riemannian metric $\tilde{g} \in C^l(\overline{Z}_2,\R^{3 \times 3})$.
We want to find a function $w \in C^{k,\alpha}(\S^2_{+})$, resp. the surface
$f(\omega) = \omega + w(\omega) \omega$, satisfying the orthogonality 
constraint
\begin{equation}
\label{eqQ0}
B[w,\tilde{g}]=\tilde{g}(\nu,\tilde{\nu}_{\R^2}) = 
\frac{1}{\sqrt{\tilde{g}^{33}}} \langle \nu, e_3 \rangle = 0,
\end{equation}
and such that $Q[w,\tilde{g}] = 0$ where $Q = Q_1,\ldots,Q_4$ is as follows:
\begin{eqnarray}
Q_1[w,\tilde{g}] & = & P^\perp W[f,\tilde{g}],\\ 
\label{eqQ2}
Q_2[w,\tilde{g}] & = & \frac{\partial H}{\partial \eta} + \tilde{h}^{\R^2}(\nu,\nu) H,\\
\label{eqQ3}
Q_3[w,\tilde{g}] & = & {\cal A}[f,\tilde{g}] - 2\pi,\\
\label{eqQ4}
Q_4[w,\tilde{g}] & = & C[f,\tilde{g}] \in \R^2. 
\end{eqnarray}
See Lemma \ref{lemmabarycenter} in the appendix for the definition 
of the twodimensional barycenter $C[f,\tilde{g}]$. We denote by 
$K= K[w,\tilde{g}]$ the space spanned by the functions
\begin{equation}
\label{eqpsifunctions}
\psi_0 = \frac{1}{\sqrt{8\pi}} H[w,\tilde{g}],\, 
\psi_i = - \sqrt{\frac{2\pi}{3}} {\rm grad}_{L^2} C^i[w,\tilde{g}] \quad (i = 1,2).
\end{equation}
A formula for $\psi_{1,2}$ is derived in (\ref{eqbary3}).
For $w = 0$, $\tilde{g} = \delta$ the functions form
an orthonormal basis of $K(0,\delta) = K_0 \subset L^2(\S^2_{+})$,
in fact by (\ref{eqbary4})
$$
\psi_0(\omega) = \frac{1}{\sqrt{2\pi}},\,
\psi_i(\omega) = \sqrt{\frac{3}{2\pi}} \langle \omega,e_i \rangle.
$$
$P^\perp$ is the $L^2$ projection, with respect to $g$, onto the orthogonal
complement of $K[w,\tilde{g}]$, thus $Q_1[w,\tilde{g}] = 0$ means
$W[f,\tilde{g}] \in K[w,\tilde{g}]$.
A function $w$ satisfying (\ref{eqQ0}) -- (\ref{eqQ4}) is called a constrained
solution for the given metric $\tilde{g}$.
The Frechet derivative $D_w Q[0,\delta]$ is the linear operator 
$$
L: C^{k,\alpha}(\S^2_{+}) \longrightarrow
C^{k-4,\alpha}(\S^2_{+}) \times C^{k-1,\alpha}(\partial \S^2_{+}) 
\times \R \times \R^2
$$
having the following components:
\begin{eqnarray}
\label{eqL1}
L_1 \varphi & = & \Delta_{\S^2} (\Delta_{\S^2} + 2) \varphi,\\
\label{eqL2}
L_2 \varphi & = &\frac{\partial}{\partial \eta} (\Delta_{\S^2} + 2)\varphi,\\
\label{eqL3}
L_3 \varphi & = & - 2 \int_{{\mathbb S}^2_{+}} \varphi\,d\mu_{\S^2},\\
\label{eqL4}
L_4^i \varphi & = & \frac{3}{2\pi} \int_{{\mathbb S}^2_{+}} 
\varphi (\omega) \langle \omega, e_i \rangle\,d\mu_{\S^2}(\omega)
\quad \mbox{ for } i = 1,2.
\end{eqnarray}
See (\ref{eqbary4}) for the derivation of (\ref{eqL4}), and note that 
$\varphi(\omega)\omega = - \varphi(\omega) \nu(\omega)$.

\begin{lemma} \label{lemmaimplicit} Let $k \geq 4$ and $l \geq k+1$. Then there 
exist open neighborhoods $W \subset C^{k,\alpha}(\S^2_{+})$ of $w \equiv 0$
and $G \subset C^{l}(Z_2,\R^{3 \times 3}_{{\rm sym}})$ of $\tilde{g} = \delta$,
and a $C^{l-k}$ function ${\bf w}:G \to W$ such that for $w \in W,\tilde{g} \in G$
\begin{equation}
\label{eqimplicit}
B[w,\tilde{g}]=0, \quad Q[w,\tilde{g}] = 0 \quad \Leftrightarrow \quad 
w = {\bf w}[\tilde{g}]. 
\end{equation}
Moreover for $\|\tilde{g} - \delta\|_{C^l(Z_2)}$ sufficiently small and 
$C = C(k,\alpha) < \infty$ we have the estimate
\begin{equation}
\label{eqimplicitestimate}
\|{\bf w}[\tilde{g}]\|_{C^{k,\alpha}(\S^2_{+})} \leq C \|\tilde{g}-\delta\|_{C^l(Z_2)}.
\end{equation}
\end{lemma}

{\em Proof. }By the coordinate expressions and the results of the appendix, 
we see that $Q[w,\tilde{g}]$ is well-defined as a map from $W \times G$ into
$C^{k-4,\alpha}(\S^2_{+}) \times C^{k-3,\alpha}(\partial \S^2_{+}) \times \R \times \R^2$,
and is of class $C^{l-k}$ under the assumptions. To construct the solution 
${\bf w}[\tilde{g}]$ we make the ansatz $w = u + \Psi[u,\tilde{g}]$, where $u \in U$, 
$\Psi[u,\tilde{g}] \in V$ are as in Lemma \ref{lemmamanifold}, in particular
$\Psi[u,\tilde{g}]$ is also of class $C^{l-k}$. The condition (\ref{eqQ0}) 
is then fulfilled, and we must solve the equation
\begin{equation}
\label{eqnonlinear}
\overline{Q}[u,\tilde{g}]: = Q[u + \Psi[u,\tilde{g}],\tilde{g}] = 0.
\end{equation}
We know from Lemma \ref{lemmamanifold} that $D_u \Psi[0,\delta] = 0$. 
Linearizing with respect to $u$ yields 
the operator $L:C^{k,\alpha}_{0}(\S^2_{+}) \to
C^{k-4,\alpha}_{\perp} \times C^{k-3,\alpha}(\partial\S^2_{+}) \times \R \times \R^2$,
where
$$
L\varphi = \left(\begin{array}{c}
\Delta_{\S^2} (\Delta_{\S^2} +2) \varphi\\
\frac{\partial (\Delta_{\S^2} \varphi)}{\partial \eta}\\
- 2 \int_{\S^2_{+}} \varphi\\
\frac{3}{2\pi} \int_{\S^2_{+}} \varphi \langle \omega, e_i \rangle\,d\mu_{\S^2}
\end{array} \right).
$$
Using Lemma \ref{lemmaiso2} it is immediate that $L$ is an 
isomorphism. By the implicit function theorem there is a solution 
$u= {\bf u}[\tilde{g}]$ of (\ref{eqnonlinear}). The $C^{l-k}$ function
${\bf w}[\tilde{g}] = {\bf u}[\tilde{g}] + \Psi\big[{\bf u}[\tilde{g}],\tilde{g}\big]$
then solves (\ref{eqimplicit}).\\
\\
Now assume that $w,\tilde{g}$ satisfy $B[w,\tilde{g}]=0$ and $Q[w,\tilde{g}] = 0$.
By uniqueness in Lemma \ref{lemmamanifold}, we then have $w = u+ \Psi[u,\tilde{g}]$ for some 
$u \in U$, and uniqueness for (\ref{eqnonlinear}) implies further $u = {\bf u}[\tilde{g}]$.
This proves the reverse implication in (\ref{eqimplicit}).\\
\\
For $\|u\|_{C^{k,\alpha}(\S^2_{+})} + \|\tilde{g}-\delta\|_{C^l(Z_2)}$ small, 
we have writing $\|\cdot\|$ for operator norms
$$
\big\|D_u \overline{Q}[u,\tilde{g}] - D_u \overline{Q}[0,\delta]\big\| 
+ \big\|D_{\tilde{g}} \overline{Q}[u,\tilde{g}] - D_{\tilde{g}} \overline{Q}[0,\delta]\big\| < \ve.
$$
We have by the fundamental theorem of calculus 
\begin{eqnarray*}
\overline{Q}[u,\tilde{g}] & = & 
D_u \overline{Q}[0,\delta] u + D_{\tilde{g}} \overline{Q}[0,\delta] (\tilde{g}-\delta)\\
&& + \int_0^1 \Big(D_u\overline{Q}[tu,(1-t)\delta + t \tilde{g}] - D_u \overline{Q}[0,\delta]\Big)\,u\,dt\\
&& + \int_0^1 \Big(D_{\tilde{g}}\overline{Q}[tu,(1-t)\delta 
+ t \tilde{g}] - D_{\tilde{g}} \overline{Q}[0,\delta]\Big)\,(\tilde{g}-\delta)\,dt.
\end{eqnarray*}
Using $\|L^{-1}\| \leq C$ and $\overline{Q}[u,\tilde{g}] = 0$ for $u = {\bf u}[\tilde{g}]$, we 
obtain after absorbing
$$
\|{\bf u}[\tilde{g}]\|_{C^{k,\alpha}(\S^2_{+})} \leq C \|\tilde{g}-\delta\|_{C^l(Z_2)}.
$$
Now $\|D_u \Psi[u,\tilde{g}\| + \|D_{\tilde{g}} \Psi[u,\tilde{g}\| \leq C$ 
for $\|u\|_{C^{k,\alpha}(\S^2_{+})} + \|\tilde{g}-\delta\|_{C^l(Z_2)}$ small, thus
$$
\|\Psi[u,\tilde{g}]\|_{C^{k,\alpha}(\S^2_+)} \leq C 
\big(\|u\|_{C^{k,\alpha}(\S^2_+)} + \|\tilde{g}-\delta\|_{C^l(Z_2)}\big).
$$
Combining yields the inequality (\ref{eqimplicitestimate}). \kasten

\begin{lemma} \label{lemmataylor} For radial graphs $f(\omega) = \omega + w(\omega) \omega$
and $l \geq 1$, consider
$$
\W:C^2(\S^2_+)\times C^l(Z_2)\rightarrow\R,\,
{\cal W}[w,\tilde{g}] = \frac{1}{4} \int_{\S^2_{+}} H^2\,d\mu_g.
$$
The functional is well-defined and of class $C^{l-1}$ on the set 
$\|w\|_{C^{1}(\S^2_+)} + \|\tilde{g}-\delta\|_{C^{0}(Z_2)} < \varepsilon_0$.
It has the derivatives, chosing $\nu(\omega) = - \omega$, 
\begin{eqnarray*}
D_w {\cal W}(0,\delta) \varphi & = & - \int_{\partial \S^2_{+}} \frac{\partial \varphi}{\partial \eta}\,ds,\\
D_{\tilde{g}}{\cal W}(0,\delta) q & = & \int_{\S^2_{+}}  \Big(-\frac{1}{2} {\rm tr}_{\S^2} q + q(\nu,\nu)
+ {\rm tr}_{\S^2} \nabla_{\!\cdot \,} q(\cdot,\nu) - \frac{1}{2} {\rm tr}_{\S^2} \nabla_\nu q \Big)\,d\mu_{\S^2}.
\end{eqnarray*}
\end{lemma}

{\em Proof. }The first formula follows from Theorem \ref{thmwillmorevariation}.
Let $\tilde{g} = \tilde{g}(\ve)$ be a family with
$$
\tilde{g}(0) = \langle\, \cdot\, ,\,\cdot\, \rangle|_{\S^2} \quad \mbox{ and } \quad
\frac{\partial \tilde{g}}{\partial \ve}|_{\ve = 0} = q.
$$
Let $\varphi:U \to \S^2_{+}$ be a parametrization. For the derivative of the normal we compute 
\begin{eqnarray*}
0 & = & \frac{\partial}{\partial \ve} \tilde{g}(\nu,\partial_\alpha \varphi)|_{\ve = 0}
= q(\nu,\partial_\alpha \varphi) + \langle \frac{\partial \nu}{\partial \ve}|_{\ve =0},\partial_\alpha \varphi\rangle,\\
0 & = &  \frac{\partial}{\partial \ve} \tilde{g}(\nu,\nu)|_{\ve = 0}
= q(\nu,\nu) + 2 \langle \frac{\partial \nu}{\partial \ve}|_{\ve =0},\nu\rangle.
\end{eqnarray*}
Thus we have
$$
\frac{\partial \nu}{\partial \ve}|_{\ve = 0} = - g^{\alpha \beta} q(\nu,\partial_\alpha \varphi) \partial_\beta \varphi
- \frac{1}{2} q(\nu,\nu) \nu.
$$
The derivative of the background connection (the Christoffel symbols) is denoted by
$\gamma(X,Y) = \frac{\partial}{\partial \ve} \tilde{D}_X Y|_{\ve = 0}$, in coordinates
$$
\gamma_{ij}^k = \frac{1}{2}(\partial_i q_{jk} + \partial_j q_{ik} - \partial_k q_{ij}).
$$
We obtain for the second fundamental form
$$
\frac{\partial h_{\alpha \beta}}{\partial \ve}|_{\ve = 0} = 
\frac{\partial}{\partial \ve} \tilde{g}(\tilde{D}_\alpha \partial_\beta \varphi,\nu)|_{\ve = 0} =
\frac{1}{2} q(\nu,\nu) h_{\alpha \beta} + \langle \gamma(\partial_\alpha \varphi,\partial_\beta \varphi),\nu\rangle.
$$
Contracting yields for the mean curvature, using $h_{\alpha \beta} = g_{\alpha \beta}$ and $H = 2$,
$$
\frac{\partial H}{\partial \ve}|_{\ve = 0} = 
- g^{\alpha \beta} q(\partial_\alpha \varphi,\partial_\beta \varphi) 
+ q(\nu,\nu) + g^{\alpha \beta}\langle \gamma(\partial_\alpha \varphi,\partial_\beta \varphi),\nu\rangle.
$$
We have further
$$
\frac{\partial}{\partial \ve} d\mu_g|_{\ve = 0 } = 
\frac{1}{2} g^{\alpha \beta} q(\partial_\alpha \varphi,\partial_\beta \varphi)\,d\mu_g.
$$
Collecting terms we find
$$
\frac{\partial}{\partial \ve}{\cal W}[0,\tilde{g}]|_{\ve = 0} = 
\int_{\S^2_{+}} \big(- \frac{1}{2} {\rm tr}_{\S^2} q + q(\nu,\nu)
+ {\rm tr}_{\S^2} \langle \gamma, \nu\rangle \big)\,d\mu_g.
$$
Finally for vectors $\tau_{1,2} \in T_\omega \S^2$ we have 
$$
\langle \gamma(\tau_1, \tau_2) ,\nu \rangle = 
\frac{1}{2} \big(\nabla_{\tau_1}q(\tau_2,\nu) + \nabla_{\tau_2}q(\tau_1,\nu) 
- \nabla_\nu q(\tau_1,\tau_2)\big).
$$
Inserting proves the second formula. \kasten

\begin{lemma}\label{lemmataylor1} For $l \geq 6$ the function ${\bf w}[\tilde{g}]$ from Lemma 
\ref{lemmaimplicit} satisfies, putting $q = \tilde{g} - \delta$,
\begin{eqnarray*}
&& \Big| {\cal W}({\bf w}[\tilde{g}],\tilde{g}) - 2\pi  
- \int_{\S^2_{+}}  \Big(-\frac{1}{2} {\rm tr}_{\S^2} q + q(\nu,\nu)
+ {\rm tr}_{\S^2} \nabla_{\!\cdot \,} q(\cdot,\nu) - \frac{1}{2} {\rm tr}_{\S^2} \nabla_\nu q \Big)\,d\mu_{\S^2}\\
&& + \int_{\partial \S^2_{+}} q(\nu,e_3)\,ds \Big|
\quad \leq \quad  C\, \|q\|^2_{C^{6}(Z_2)}. 
\end{eqnarray*}
\end{lemma}

{\em Proof. }Putting $\varphi = D{\bf w}[\,\delta\,] q$ we compute
\begin{eqnarray*}
&& D_w {\cal W}[0,\delta] \varphi + D_{\tilde{g}} {\cal W}[0,\delta] q\\
& = & \int_{\S^2_{+}}  \Big(-\frac{1}{2} {\rm tr}_{\S^2} q + q(\nu,\nu)
+ {\rm tr}_{\S^2} \nabla_{\!\cdot \,} q(\cdot,\nu) - \frac{1}{2} {\rm tr}_{\S^2} \nabla_\nu q \Big)\,d\mu_{\S^2}
- \int_{\S^2_{+}} \frac{\partial \varphi}{\partial \eta}\,ds.
\end{eqnarray*}
On the other hand we had in Lemma \ref{lemmamanifold} 
$$
0 = D_w B[0,\delta] \varphi + D_{\tilde{g}} B[0,\delta] q
= \frac{\partial \varphi}{\partial \eta} - q(\nu,e_3).
$$
The claim follows by Taylor's formula, taking $k = 4, l= 6$ in Lemma \ref{lemmaimplicit}.
\kasten

\section{Blowup at boundary points}
Let $\Omega \subset \R^3$ be a bounded domain of class $C^m$, $m \geq 7$, 
with boundary $S = \partial \Omega$. At a given point $a \in S$ we let $N(a)$ 
be the interior unit normal and choose an orthonormal basis $v_1(a),v_2(a)$ of $T_aS$. 
For $r_0 = r_0(\Omega) > 0$ we have a graph representation
\begin{equation}
\label{graph1}
f^a:D_{r_0} \to \R^3,\, f^a(x) = a + x^1 v_1(a) + x^2 v_2(a) + \varphi^a(x) N(a), 
\end{equation}
such that 
\begin{equation}
\label{graphbounds}
\|\varphi^a\|_{C^{m}(D_{r_0})} \leq C = C(\Omega).
\end{equation}
Since $\varphi^a(0) = 0$ and $D\varphi^a(0) = 0$ we have 
\begin{equation}
\label{graphgrowth}
|\varphi^a(x)| \leq C |x|^2 \quad \mbox{ and } \quad  |D\varphi^a(x)| \leq C |x|. 
\end{equation}
We extend the graph parametrization to a diffeomorphism 
$$
F^a:Z_{r_0}= D_{r_0} \times (-r_0,r_0) \to \R^3,\, F^a(x,z) = 
f^a(x) + z N(a).
$$
Using indices $i,j,k = 1,2$ we compute for $\tilde{g}^a = (F^a)^\ast \langle \,\cdot\,,\,\cdot\, \rangle_{\R^3}$
\begin{eqnarray}
\nonumber
\partial_i F^a(x,z) & = & v_i(a) + \partial_i \varphi^a(x)N(a)\\ 
\nonumber
\partial_3 F^a(x,z) & = & N(a),\\
\tilde{g}^a_{ij}(x,z)     & = & \delta_{ij} + \partial_i \varphi^a(x) \partial_j \varphi^a(x)\\
\nonumber
\tilde{g}^a_{i3}(x,z)     & = & \partial_i \varphi^a(x)\\ 
\nonumber
\tilde{g}^a_{33}(x,z) & = & 1.
\end{eqnarray}
Next consider the dilations 
$$
\sigma_\lambda:\R^3 \to \R^3,\, \sigma_\lambda(x,z) = 
(\lambda x,\lambda z) \mbox{ where } \lambda > 0.
$$
Clearly $\sigma_\lambda(Z_R) \subset Z_{r_0}$ for $R \leq \frac{r_0}{\lambda}$.
We obtain the Riemannian isometry
\begin{equation}
\label{eqchart}
F^{a,\lambda}:\big(Z_{\frac{r_0}{\lambda}},\lambda^2 \tilde{g}^{a,\lambda}\big) \to 
F^a(Z_{r_0}) \subset \R^3,\,
F^{a,\lambda}(x,z) = F^a(\lambda x,\lambda z),
\end{equation}
where the metric $\tilde{g}^{a,\lambda}$ is given by 
\begin{equation}
\label{eqpullback}
\tilde{g}^{a,\lambda}:Z_{\frac{r_0}{\lambda}} \to \R^{3 \times 3}, 
\tilde{g}^{a,\lambda}(x,z) = \lambda^{-2} (\sigma_\lambda)^\ast \tilde{g}^a(x,z)
= \tilde{g}^a(\lambda x,\lambda z).
\end{equation}
The metric satisfies, as a function of $(\lambda,x,z)$ for $a \in S$ fixed, 
$$
\tilde{g}^{a,\lambda}(x,z) \in C^{m-1}(\big[0,\frac{r_0}{2}\big] \times Z_2,\R^{3 \times 3}) \quad 
\mbox{ where } \tilde{g}^{a,0}_{ij} = \delta_{ij}.
$$
Moreover the above expansions yield bounds, for a constant $C = C(\Omega)$, 
\begin{eqnarray}
\nonumber
\|\tilde{g}^{a,\lambda}_{ij} - \delta_{ij}\|_{C^{m-1}(Z_2)}
& \leq & C \lambda^2,\\
\label{eqmetricasymptotics}
\|\tilde{g}^{a,\lambda}_{i3}\|_{C^{m-1}(Z_2)} & \leq & C \lambda,\\
\nonumber
\tilde{g}^{a,\lambda}_{33} - 1 & \equiv  & 0.
\end{eqnarray}
We compute more precisely
\begin{equation}
\label{eqmetricexpansion1} q_{ij}(x,z):=
\frac{\partial}{\partial \lambda}\tilde{g}^{a,\lambda}_{ij}(x,z)|_{\lambda = 0} = 
\begin{cases}
0 & \mbox{ for } 1 \leq i,j \leq 2\\
h^S_{ik}(a)x^k & \mbox{ for } i = 1,2 \mbox{ and } j=3\\
0 & \mbox{ for } i = j = 3.
\end{cases}
\end{equation}
Taylor expansion yields for $C = C(\Omega)$
\begin{equation}
\label{eqtaylor1}
\|\tilde{g}^{a,\lambda} - (\delta_{ij} + \lambda q_{ij})\|_{C^{m-3}(Z_2)} 
\leq C \lambda^2 \quad \mbox{ where }0 \leq \lambda \leq \lambda_0(\Omega).
\end{equation}

\begin{lemma} \label{lemmaintegrals1} For $q_{ij}(x,z)$ as in (\ref{eqmetricexpansion1}), 
we have the following formulae:
\begin{eqnarray*}
\int_{\S^2_{+}} q(\nu,\nu)\,d\omega & = & \frac{\pi}{2} H^S(a),\\
\int_{\S^2_{+}} {\rm tr}_{\S^2} q\,d\omega & = & - \frac{\pi}{2} H^S(a),\\
\int_{\S^2_{+}} {\rm tr}_{\S^2} \nabla_\nu q\,d\omega & = & \frac{\pi}{2} H^S(a),\\ 
\int_{\S^2_{+}} {\rm tr}_{\S^2}\nabla_{\cdot} q(\cdot,\nu)\,d\omega & = & 
- \frac{\pi}{2} H^S(a),\\
\int_{\partial \S^2_{+}} q(\nu,e_3)\,ds & = & - \pi H^S(a).
\end{eqnarray*}
\end{lemma}

{\em Proof. }We compute writing $\omega = (\sin \theta\, \xi,\cos \theta)$ 
for $\xi \in \S^1$, $0 \leq \theta \leq \frac{\pi}{2}$ 
$$
\int_{\S^2_{+}} q(\nu,\nu)\,d\omega = 
2 \int_0^{\frac{\pi}{2}} \int_{\S^1} h^S(a)(\xi,\xi) \sin^3 \theta \cos \theta \,d\xi d\theta
= \frac{\pi}{2} H^S(a).
$$
Since ${\rm tr\,}q = 0$ we get
$$
\int_{\S^2_{+}} {\rm tr}_{\S^2} q\,d\omega = 
- \int_{\S^2_{+}} q(\nu,\nu)\,d\omega = - \frac{\pi}{2} H^S(a).
$$
Differentiating the equation $q(t\omega) = t q(\omega)$ at $t = 1$, we get 
$$
\int_{\S^2_{+}} {\rm tr}_{\S^2} \nabla_\nu q\,d\omega = 
- \int_{\S^2_{+}} {\rm tr}_{\S^2} q\,d\omega =  \frac{\pi}{2} H^S(a).
$$
Now we compute using the definition of $q$
$$
{\rm tr}_{\S^2} \nabla_{\cdot} q(\cdot,\nu) = 
{\rm tr}_{\R^3} \nabla_{\cdot} q(\cdot,\nu) - \nabla_\nu q(\nu,\nu)
= \langle \nu,e_3 \rangle \sum_{i=1}^2 \partial_i q_{i3} + q(\nu,\nu).
$$
Using $\sum_{i=1}^2 \partial_i q_{i3} = h^S_{11}(a) + h^S_{22}(a) = H^S(a)$ we see
$$
\int_{\S^2_{+}} {\rm tr}_{\S^2} \nabla_{\cdot} q(\cdot,\nu)\,d\omega = 
- \int_0^{\frac{\pi}{2}} \int_{\S^1} H^S(a) \cos \theta \sin \theta\,d\xi d\theta 
+ \frac{\pi}{2} H^S(a) = - \frac{\pi}{2} H^S(a).
$$
Finally we compute the boundary integral
$$
\int_{\partial \S^2_{+}} q(\nu,e_3)\,ds =
- \int_{\S^1} h(a)(\xi,\xi)\,d\xi = - \pi H(a).
$$
\kasten

For $l = m-1 \geq 6$ and $0 \leq \lambda \leq \lambda_0(\Omega)$ the metric 
$\tilde{g}^{a,\lambda}$ belongs to the neighborhood $G$ of the standard 
metric as in Lemma \ref{lemmaimplicit}.  We put 
$w^{a,\lambda} = {\bf w}[\tilde{g}^{a,\lambda}]$ and
$q^{a,\lambda} = \tilde{g}^{a,\lambda} - \delta$. The Taylor 
expansion from Lemma \ref{lemmataylor1} then yields 
\begin{eqnarray*}
\label{eqwillmoreexpansion1}
\Big|{\cal W}\big(w^{a,\lambda},\tilde{g}^{a,\lambda}\big) - 2\pi & - &
\int_{\S^2_{+}} \Big(-\frac{1}{2} {\rm tr}_{\S^2} q^{a,\lambda} + q^{a,\lambda}(\nu,\nu)
+ {\rm tr}_{\S^2} \nabla_{\cdot} q(\cdot,\nu) - \frac{1}{2} {\rm tr}_{\S^2} \nabla_\nu q\Big)\,d\omega\\
& + & \int_{\partial \S^2_{+}} q^{a,\lambda}(\nu,e_3)\,ds \Big| \leq 
C\,\Vert q^{a,\lambda}\Vert^2_{C^{6}(Z_2)}.
\end{eqnarray*}
Now $\|q^{a,\lambda}-\lambda q\|_{C^{m-2,\alpha}(Z_2)} \leq C \lambda^2$ by (\ref{eqtaylor1}),
hence evaluating the integrals shows 
\begin{equation}
\label{eqwillmoreexpansion1}
\Big|\frac{{\cal W}\big(w^{a,\lambda},\tilde{g}^{a,\lambda}\big) - 2\pi}{\lambda} + \pi H(a) \Big| 
\leq C \lambda \quad \mbox{ where }C = C(\Omega).
\end{equation}
Transforming back yields the following result, where by ${\cal M}^{k,\alpha}(S)$ 
we denote the set of $C^{k,\alpha}$ immersions of $\S^2_{+}$ meeting $S$ orthogonally 
from inside along the boundary.

\begin{proposition} \label{propexpansion} Let $\Omega$ be of class $C^m$ for 
$m \geq 7$, and $k: = m-2$. Then for $a \in S$ and $0 < \lambda \leq \lambda_0$ 
the $\phi^{a,\lambda}(\omega)  = F^{a,\lambda}\big((1+w^{a,\lambda}(\omega))\omega\big)$ 
belong to ${\cal M}^{k,\alpha}(S)$, have area ${\cal A}(\phi^{a,\lambda}) = 2\pi \lambda^2$, 
are centered at $a \in S$ and satisfy 
$$
\Big|\frac{{\cal W}(\phi^{a,\lambda})- 2\pi}{\lambda} + \pi H(a) \Big| \leq C \lambda
\quad \mbox{ where }C = C(\Omega).
$$
In particular we see that $\inf_{f \in {\cal M}(S)} {\cal W}(f) < 2\pi$.
\end{proposition}

{\bf Remark. }Suppose that a sequence of immersions $f_k \in {\cal M}(S)$ satisfies 
\begin{equation}
\label{eqassumedbounds}
{\rm diam\,}f_k(D) \to 0,\quad 
{\cal A}(f_k) \leq C, \quad  {\rm L}(f_k|_{\partial \Sigma}) \leq C.
\end{equation}
It is not difficult to show that then $\liminf_{k \to \infty} {\cal W}(f_k) \geq 2\pi$.
Thus for a ${\cal W}$-minimizing sequence $f_k$ in ${\cal M}(S)$ one of the bounds 
in (\ref{eqassumedbounds}) must be violated in view of Proposition \ref{propexpansion}.
For $\Omega$ convex the length bound could in fact be dropped using the Gau{\ss} 
Bonnet theorem. Global bounds for the Willmore energy of surfaces with free boundary 
are proved in recent work by Volkmann \cite{Vol14}.\\
\\
In the following lemma we check how the constrained solutions transform 
when changing the orthonormal basis $v_1(a),v_2(a)$ used to identify 
$T_a S$ with $\R^2$. 

\begin{lemma} Let $w^{a,\lambda}$ be the solution with respect
to the basis $v_{1,2} = v_{1,2}(a)$ of $T_aS$, and let $T \in {\mathbb SO}(2)$.
Then the corresponding solution $w^{T,a,\lambda}$ with respect to
the basis $v^T_j = T_{ij} v_i$ is given by 
$$
w^{T,a,\lambda} = w^{a,\lambda} \circ T, \quad \mbox{ where we identify }\,
T \hat{=} \left(\begin{array}{cc} T & 0\\ 0 & 1 \end{array} \right).
$$
In particular we have $\phi^{T,a,\lambda} = \phi^{a,\lambda} \circ T$.
\end{lemma}

{\em Proof. }We compute
\begin{eqnarray*}
f^a(Tx) & = &a + (Tx)^1 v_1 + (Tx)^2 v_2 + \varphi(Tx)N(a)\\
& = & a + (T_{11} x^1 + T_{12} x^2) v_1 + (T_{21}x^1 + T_{22} x^2) v_2 + \varphi^a(Tx) N(a)\\
& = & a + x^1 (\underbrace{T_{11} v_1 + T_{21} v_2}_{= v_1^T}) 
+ x^2 (\underbrace{T_{12} v_1 + T_{22} v_2}_{=v_2^T}) + \varphi^a(Tx) N(a).
\end{eqnarray*}
This shows $\varphi^{T,a}(x) = \varphi^a(Tx)$ and $f^{T,a}(x) = f^a(Tx)$ on $D_{r_0}$. 
It follows that $F^{T,a,\lambda}(x,z) = F^{a,\lambda}(Tx,z)$ and hence 
$$
\tilde{g}^{T,a,\lambda} 
= \lambda^{-2} (F^{T,a,\lambda})^\ast \langle \,\cdot \,,\,\cdot\,\rangle
= \lambda^{-2} (F^{a,\lambda} \circ T)^\ast \langle \,\cdot \,,\,\cdot\,\rangle
= T^\ast \tilde{g}^{a,\lambda}.
$$
The boundary value problem (\ref{eqimplicit}) is Riemannian 
invariant, that is
\begin{eqnarray*}
B[w^{a,\lambda} \circ T,\tilde{g}^{T,a,\lambda}] & = & B[w^{a,\lambda},\tilde{g}^{a,\lambda}] \circ T = 0,\\
Q^i[w^{a,\lambda} \circ T,\tilde{g}^{T,a,\lambda}] & = & Q^i[w^{a,\lambda},\tilde{g}^{a,\lambda}] \circ T = 0,
\quad \mbox{ for } i = 1,2,\\
Q^3[w^{a,\lambda} \circ T,\tilde{g}^{T,a,\lambda}] & = & Q^3[w^{a,\lambda},\tilde{g}^{a,\lambda}] = 0,\\
Q^4[w^{a,\lambda} \circ T,\tilde{g}^{T,a,\lambda}] & = & T^{-1} Q^4[w^{a,\lambda},\tilde{g}^{a,\lambda}] = 0.
\end{eqnarray*}
By uniqueness in Lemma \ref{lemmaimplicit} we conclude that $w^{T,a,\lambda} = w^{a,\lambda} \circ T$.
\kasten

We now study the reduced energy function 
\begin{equation}
\label{eqreducedenergy} 
\Wbar:S \times [0,\lambda_0) \to \R,\,\Wbar(a,\lambda) 
= {\cal W}(\phi^{a,\lambda}) 
= {\cal W}(w^{a,\lambda},\tilde{g}^{a,\lambda}).
\end{equation}
We already know that 
$$
\Wbar(a,0) \equiv 2\pi,\, \nabla_a \Wbar(a,0) \equiv 0 \quad \mbox{and}\quad 
\frac{\partial \Wbar}{\partial \lambda}(a,0) = - \pi H(a).
$$ 
For further computations we assume w.l.o.g. that $0 \in S$, $N(0) = e_3$,
and chose an orthonormal frame $v_{1,2} \in C^{m-1}(U,\R^3)$
on a neighborhood $U \subset S$ such that
\begin{equation}
\label{eqparallel}
D_{v_i} v_j(0) = h^S(0)(v_i,v_j) N^S(0). 
\end{equation}
The $v_i(a)$ can be obtained for instance by Gram-Schmidt applied to the 
coordinate vectors of the local graph representation. In order to have 
$\Wbar$ of class $C^r$ for $r \geq 1$, we assume in the following that 
$m = 6 +2r$. Taking $k = 4$, one then checks that the map 
$$
S \times [0,\lambda_0] \to C^{4+r}(Z_2,\R^{3 \times 3}),\,(a,\lambda) \mapsto \tilde{g}^{a,\lambda},
$$
is of class $C^r$, which implies also $\Wbar \in C^r(S \times [0,\lambda_0])$.

For example, for $m = 10$ we can take $r = 2$ and deduce 
\begin{equation}
\label{eqmixedderivatives}
\frac{\partial}{\partial \lambda} \nabla_a \Wbar(a,0) = 
\nabla_a \frac{\partial \Wbar}{\partial \lambda} (a,0) = 
- \pi \nabla H(a). 
\end{equation}

\begin{theorem} \label{thmreduction} Let $\Omega \subset \R^3$ be a bounded 
domain of class $C^8$. Put $S = \partial \Omega$ and
$$
{\cal M}^{4,\alpha}_{\lambda}(S) = \{f \in {\cal M}^{4,\alpha}(S): {\cal A}(f) = 2\pi \lambda^2\}.
$$
For any $\lambda \in (0,\lambda_0]$ and $a \in S$ the following are equivalent:
\begin{itemize}
\item[{\rm (1)}] $a$ is a critical point of $\Wbar(\cdot,\lambda)$.
\item[{\rm (2)}] $\phi^{a,\lambda}$ is a critical point of the Willmore 
functional in ${\cal M}^{4,\alpha}_{\lambda}(S)$.
\item[{\rm (3)}] $\phi^{a,\lambda}$ solves the boundary value problem
\begin{eqnarray*}
\Delta H + |A^\circ|^2 H = & \alpha H &\mbox{ for some } \alpha \in \R,\\
\frac{\partial  f}{\partial \eta} = & N^S \circ f & \mbox{ along }\partial \S^2_{+},\\
\frac{\partial H}{\partial \eta} + h^S(\nu,\nu) H = & 0 & \mbox{ along }\partial \S^2_{+}.
\end{eqnarray*}
\end{itemize}
\end{theorem}

\begin{corollary} \label{corexistence} Let $\Omega \subset \R^3$
be a bounded domain of class $C^8$. Then for any $\lambda \in (0,\lambda_0]$ 
there exist two different critical points of the Willmore functional in
${\cal M}^{4,\alpha}_{\lambda}(S)$, corresponding to the extrema of the 
function $\Wbar(\cdot,\lambda)$.
\end{corollary}

{\em Proof. } From Proposition \ref{propexpansion} we have for $a_1,a_2 \in S$ 
$$
\pi |H^S(a_1) - H^S(a_2)| \leq 
\frac{|\Wbar(a_1,\lambda) - \Wbar(a_2,\lambda)|}{\lambda} + C \lambda.
$$
If there is a sequence $\lambda_k \searrow 0$ such that each function
$\Wbar(\cdot,\lambda_k)$ is constant, then $H^S$ must be constant
and hence $\Omega$ is a round ball by Alexandroffs theorem. By
symmetry we then have infinitely many critical points. On the other hand,
if $\Wbar(\cdot,\lambda)$ is not constant, then it attains its extrema
at different points $a_1(\lambda),a_2(\lambda) \in S$. The surfaces
$\phi^{a_i(\lambda),\lambda}$ are then geometrically different, since
the $a_i(\lambda)$ are their barycenters. \kasten

As noted in \cite{PX09,Tak68} the number of critical points is in fact bounded 
below by the Ljusternik-Shnirelman category of $S$, which equals three 
if $S$ is a surface of higher genus. We have also the following fact 
about the concentration points for $\lambda \searrow 0$.

\begin{corollary} \label{corconcentrationpoints}  Let $\Omega \subset \R^3$ 
be a bounded domain of class $C^{10}$, and assume that the $\phi^{a_k,\lambda_k}$
are critical points of the Willmore functional in ${\cal M}_{\lambda_k}(S)$,
where $\lambda_k \to 0$ and $a_k \to a \in S$. Then $\nabla H^S(a) = 0$.
\end{corollary}

{\em Proof. }We have $\nabla_a \Wbar(a_k,\lambda_k) = 0$ 
by assumption. Using $\nabla_a \Wbar(a,0) \equiv 0$ which follows from 
$\Wbar(a,0) \equiv 2\pi$, we get
$$
0 = \frac{\nabla_a \Wbar(a_k,\lambda_k) - \nabla_a \Wbar(a_k,0)}{\lambda_k}
= \mint_0^{\lambda_k} \frac{\partial}{\partial \lambda} \nabla_a \Wbar(a_k,\lambda)\,d\lambda
\to \frac{\partial}{\partial \lambda} \nabla_a \Wbar(a,0). 
$$
Claim (1) follows from (\ref{eqmixedderivatives}). \kasten

{\em Proof of Theorem \ref{thmreduction}. }For $\lambda_0 > 0$ sufficiently small, we show that critical 
points of $\Wbar(\cdot,\lambda)$, $\lambda \in (0,\lambda_0]$, correspond to critical 
points of the Willmore functional in ${\cal M}^{4,\alpha}_\lambda(S)$.
Consider the constrained solutions 
$$
\phi^{a,\lambda}(\omega) = F^{a,\lambda}\big(\omega + w^{a,\lambda}(\omega)\omega\big),\,
\quad \omega \in \S^2_{+}.
$$
For fixed $\lambda$ the family $\phi^{a,\lambda}$ is a variation in 
${\cal M}^{4,\alpha}_\lambda(S)$. Now 
\begin{eqnarray*}
H[\phi^{a,\lambda}] & = &
H[w^{a,\lambda},(F^{a,\lambda})^\ast \langle \cdot,\cdot \rangle] 
= \lambda^{-1} H[w^{a,\lambda},\tilde{g}^{a,\lambda}],\\
W[\phi^{a,\lambda}] & = & W[w^{a,\lambda},(F^{a,\lambda})^\ast \langle \cdot,\cdot \rangle]
= \lambda^{-3} W[w^{a,\lambda},\tilde{g}^{a,\lambda}]. 
\end{eqnarray*}
Thus for $\psi_i = \psi_i[w^{a,\lambda}, \tilde{g}^{a,\lambda}]$, $i = 0,1,2$, as in 
(\ref{eqpsifunctions}) we have 
$$
W[\phi^{a,\lambda}] \in {\rm Span\,}\{\psi_0,\psi_1,\psi_2\}. 
$$
We further have along $\partial \S^2_{+}$
$$
\langle \nu[\phi^{a,\lambda}],N^S \circ \phi^{a,\lambda} \rangle = 0 \quad \mbox{ and } 
\quad \Big(\frac{\partial H}{\partial \eta} + h^S(\nu,\nu) H\Big)[\phi^{a,\lambda}] = 0.
$$
For the two-dimensional barycenter defined in Lemma \ref{lemmabarycenter} we see 
$$
\pi_S \Big(\mint_{\S^2_{+}} \phi^{a,\lambda}\,d\mu_{\phi^{a,\lambda}}\Big)
= F^{a,\lambda}\big(C[w^{a,\lambda},\tilde{g}^{a,\lambda}]\big)
= F^{a,\lambda}(0) = a.
$$
This summarizes the conditions for constrained solutions. Next we 
study variations corresponding to the parameter $a$.\\
\\
Assume that $0 \in S$, $N^S(0) = e_3$, is a critical point for the function 
$\Wbar^\lambda = \Wbar(\cdot,\lambda)$. 
Choose an orthonormal frame $v_{1,2}(a) \in C^{7}(U,\R^3)$ nearby, 
such that $\nabla^S_{v_i} v_j(0) = (D_{v_i} v_j)^\top(0) = 0$.  
The map $F^{a,\lambda}$ is given explicitely by
$$
F^{a,\lambda}(x,z) = a + \lambda x^a + (\varphi^a(\lambda x) + \lambda z)N^S(a),
\quad \mbox{ where } x^a = x^1 v_1(a) + x^2 v_2(a). 
$$
Taking the derivative $\frac{\partial}{\partial a^i}$ at $a = 0$ gives
$$
\frac{\partial F^{a,\lambda}}{\partial a^i}(x,z)|_{a=0} = 
e_i + (\varphi^0(\lambda x) + \lambda z) W^S(0) e_i
+ \lambda \frac{\partial x^a}{\partial a^i}(x)|_{a=0} 
+ \frac{\partial \varphi^a}{\partial a^i}(\lambda x)|_{a=0} e_3.
$$
We have $\frac{\partial v_j}{\partial a_i}(0) = h^S_{ij}(0) e_3$, thus
$$
\frac{\partial x^a}{\partial a^i}(x)|_{a=0} = \big(x^1 h^S_{1i}(0) + x^2 h^S_{2i}(0)\big) e_3.
$$
Next we write $\frac{\partial \varphi^a}{\partial a^i}|_{a=0}$ 
in terms of the graph function $\varphi^{a=0}$, using the equation
$$
\langle F^{a,\lambda}(x,0),e_3 \rangle = \varphi^0\big(\pi_{\R^2} F^{a,\lambda}(x,0)\big).
$$
The derivative $\frac{\partial}{\partial a^i}$ yields at $a= 0$
$$
\lambda \big(x^1 h^S_{1i}(0) + x^2 h^S_{2i}(0)\big) 
+ \frac{\partial \varphi^a}{\partial a_i}|_{a=0} =
\big \langle \nabla \varphi^0(\lambda x), e_i + \varphi^0(\lambda x) W^S(0)e_i \big \rangle. 
$$
Rearranging gives 
\begin{equation}
\label{eqgraphparameterderivative}
\frac{\partial \varphi^a}{\partial a_i}|_{a=0} = 
- \lambda (x^1 h^S_{1i}(0) + x^2 h^S_{2i}(0))
+ (\delta_{ij} - \varphi^0(\lambda x) h^S_{ij}(0)) \partial_j \varphi^0(\lambda x).
\end{equation}
Reinserting yields the formula
\begin{eqnarray}
\label{eqpseudotranslation}
\frac{\partial F^{a,\lambda}}{\partial a^i}(x,z)|_{a=0} & = & 
e_i - \big(\varphi^0(\lambda x) + \lambda z\big) h^S_{ij} (0) e_j\\
\nonumber
&& + \big(\delta_{ij} - \varphi^0(\lambda x) h^S_{ij}(0)\big) \partial_j \varphi^0(\lambda x) e_3.
\end{eqnarray}
By the assumptions on $\varphi^0$ we have 
$$
\frac{1}{\lambda} \|\varphi^0(\lambda x)\|_{C^8(B_2)} 
+ \|D\varphi^0(\lambda x)\|_{C^{7}(B_2)} 
\leq C \lambda \quad 
\mbox{ for } \lambda \leq \frac{r_0}{2}.
$$
This implies
$$
\Big\|\frac{\partial F^{a,\lambda}}{\partial a^i}|_{a=0} - e_i\Big\|_{C^{7}(Z_2)} \leq C \lambda.
$$
Now consider the $\phi^{a,\lambda} = F^{a,\lambda} \circ f^{a,\lambda}$ where 
$f^{a,\lambda}(\omega) = \big(1+ w^{a,\lambda}(\omega)\big) \omega$. We have
$$
\frac{\partial \Wbar}{\partial a_i}(0,\lambda) 
= \frac{\partial}{\partial a_i}{\cal W}(\phi^{a,\lambda})|_{a = 0}
= D{\cal W}(\phi^{0,\lambda}) \cdot \Big( 
\underbrace{\frac{\partial F^{a,\lambda}}{\partial a^i}|_{a=0} \circ f^{0,\lambda}
+ DF^{0,\lambda} \circ f^{0,\lambda} \frac{\partial f^{a,\lambda}}{\partial a^i}|_{a=0}}_{=: Y_i}\Big).
$$
We transform back to the reference chart, defining the vector field
$$
X_i:\S^2_{+} \to \R^3,\,X_i(\omega) = 
\lambda DF^{0,\lambda}\big(f^{0,\lambda}(\omega)\big)^{-1} \cdot Y_i (\omega).
$$
By Riemannian invariance, we then get 
$$
D_w{\cal W}(f^{0,\lambda},\tilde{g}^{0,\lambda}) \cdot X_i =
\lambda\, D{\cal W}(\phi^{0,\lambda}) \cdot Y_i.
$$
We want to show that $X_i \approx e_i$ for sufficiently small $\lambda > 0$. From 
the definition $F^0(x,z) = (x,z+\varphi^0(x))$ we see $(F^0)^{-1}(x,z) = (x,z-\varphi^0(x))$, thus
$$
(DF^0)^{-1}(\lambda x, \lambda z) 
= D\big((F^0)^{-1}\big)(F^0(\lambda x, \lambda z)) 
= {\rm Id} - d\varphi^0(\lambda x) \otimes e_3.
$$
Now $\lambda DF^{0,\lambda}(x,z) = DF^0(\lambda x,\lambda z)$,  which yields 
\begin{eqnarray*}
\lambda (DF^{0,\lambda})^{-1} \frac{\partial F^{a,\lambda}}{\partial a^i}|_{a=0} & = &
e_i - \big(\varphi^0(\lambda x) + \lambda z\big) h^S_{ij} (0) e_j\\
\nonumber
&& + \big(\delta_{ij} - \varphi^0(\lambda x) h^S_{ij}(0)\big) \partial_j \varphi^0(\lambda x) e_3\\
&& - \big(\partial_i \varphi^0(\lambda x) 
- \big(\varphi^0(\lambda x) + \lambda z\big) h^S_{ij} (0) \partial_j\varphi^0(\lambda x)\big) e_3\\
& = & e_i - \big(\varphi^0(\lambda x) + \lambda z\big) h^S_{ij} (0) e_j 
+ \lambda z h^S_{ij}(0) \partial_j \varphi^0(\lambda x) e_3.
\end{eqnarray*}
In particular
\begin{equation}
\label{eqtranslation1}
\Big\|\lambda (DF^{0,\lambda})^{-1} \frac{\partial F^{a,\lambda}}{\partial a^i}|_{a=0} 
- e_i\Big\|_{C^{7}(Z_2)}
\leq C \lambda.
\end{equation}
The functions $w^{a,\lambda}$ are defined as the solutions of the equation 
$Q[w,\tilde{g}^{a,\lambda}] = 0$, taking $k = 4$ in Lemma \ref{lemmaimplicit}.
From (\ref{eqimplicitestimate}) we have the bound
$$
\|w^{0,\lambda}\|_{C^{4,\alpha}(\S^2_{+})} 
\leq C \|\tilde{g}^{0,\lambda}-\delta\|_{C^{5}(Z_2)}
\leq C \lambda.
$$
To estimate $\frac{\partial w^{a,\lambda}}{\partial a^i}|_{a=0}$, we compute
$$
0 = \frac{\partial}{\partial a^i} Q[w^{a,\lambda},\tilde{g}^{a,\lambda}]|_{a=0} 
= D_w Q[w^{0,\lambda},\tilde{g}^{0,\lambda}] 
\cdot \frac{\partial w^{a,\lambda}}{\partial a^i}|_{a=0}
+ D_{\tilde{g}}Q[w^{0,\lambda},\tilde{g}^{0,\lambda}] \cdot
\frac{\partial \tilde{g}^{a,\lambda}}{\partial a^i}|_{a=0}. 
$$
For $\lambda > 0$ sufficiently small we have $\|\tilde{g}^{0,\lambda}-\delta\|_{C^5(Z_2)}$ 
small and hence $\|w^{0,\lambda}\|_{C^{4,\alpha}(\S^2_{+})}$ small, so that 
$D_w Q[w^{0,\lambda},\tilde{g}^{0,\lambda}]$ is close to the invertible 
Operator $L = D_w Q[0,\delta]$. Thus we can estimate 
\begin{eqnarray*}
\big\|\frac{\partial w^{a,\lambda}}{\partial a^i}|_{a=0}\big\|_{C^{4,\alpha}(\S^2_{+})}
& \leq & C\, \big\|D_{\tilde{g}}Q[w^{0,\lambda},\tilde{g}^{0,\lambda}] \cdot
\frac{\partial \tilde{g}^{a,\lambda}}{\partial a^i}|_{a=0}\|_{
C^{0,\alpha}(\S^2_{+}) \times C^{1,\alpha}(\partial \S^2_{+}) \times \R \times \R^2}\\
& \leq & C\,\|\frac{\partial \tilde{g}^{a,\lambda}}{\partial a^i}|_{a=0}\|_{C^{5}(Z_2)}\\
& \leq & C \lambda.
\end{eqnarray*}
In the last estimate we used the definition of $\tilde{g}^{a,\lambda}$, the 
formula (\ref{eqgraphparameterderivative}) and the $C^8$ bound on $\varphi^0$.
For $f^{a,\lambda}(\omega) = \big(1+ w^{a,\lambda}(\omega)\big) \omega$ we obtain 
\begin{equation}
\label{eqtranslation2}
\|f^{0,\lambda}(\omega) - \omega\|_{C^{4,\alpha}(\S^2_{+})} +
\|\frac{\partial f^{a,\lambda}}{\partial a^i}|_{a=0}\|_{C^{4,\alpha}(\S^2_{+})}
\leq C \lambda.
\end{equation}
Now we have 
$$
X_i(\omega) = 
\Big(\lambda (DF^{0,\lambda})^{-1} \frac{\partial F^{a,\lambda}}{\partial a^i}|_{a=0}\Big)|_{f^{0,\lambda}(\omega)}
+ \lambda \frac{\partial f^{a,\lambda}}{\partial a^i}|_{a=0}.
$$
Combining (\ref{eqtranslation1}) and (\ref{eqtranslation2}) we conclude
\begin{equation}
\label{eqtranslation3}
\|X_i - e_i\|_{C^{4,\alpha}(\S^2_{+})} \leq C \lambda.
\end{equation}
Now write $Y_i = \varphi_i \nu_{\phi^{0,\lambda}} + D\phi^{0,\lambda} \tau_i$, and compute
$$
\frac{\partial}{\partial a_i} {\cal W}(\phi^{a,\lambda})|_{a=0} = 
\frac{1}{2} \int_{\S^2_{+}} W[\phi^{0,\lambda}] \varphi_i\,d\mu_g
+ \frac{1}{2} \int_{\partial \S^2_{+}} \Big(\varphi_i \frac{\partial H}{\partial \eta}
- \frac{\partial \varphi_i}{\partial \eta} H - \frac{1}{2} g(\tau_i,\eta)\Big)\,ds_g.
$$
As $\phi^{a,\lambda}(\partial \S^2_{+}) \subset S$ the vector $Y_i(\omega)$ is tangent
to $S$ at $\phi^{0,\lambda}(\omega)$, for any $\omega \in \partial \S^2_{+}$. Since 
$\frac{\partial \phi^{0,\lambda}}{\partial \eta} = N^S \circ \phi^{0,\lambda}$
we get $g(\tau_i,\eta) \equiv 0$ along $\partial \S^2_{+}$. Furthermore
\begin{eqnarray*}
\frac{\partial H}{\partial \eta} & = & - h^S(\nu,\nu) H \quad \mbox{(boundary condition for $\phi^{a,\lambda}$)}\\
\frac{\partial \varphi_i}{\partial \eta} & = & - h^S(\nu,\nu) \varphi_i \quad 
\mbox{ (admissibility as in (\ref{eqlinboundary4}))}.
\end{eqnarray*}
Thus all boundary terms cancel and we get putting $\xi_i = \tilde{g}^{0,\lambda}\big(X_i,\nu_{f^{0,\lambda}}\big)$
\begin{eqnarray*}
\frac{\partial}{\partial a_i} {\cal W}(\phi^{a,\lambda})|_{a=0} & = &
\frac{1}{2} \int_{\S^2_{+}} \langle \vec{W}[\phi^{0,\lambda}],Y_i \rangle \,d\mu_{\phi^{0,\lambda}}\\
& = & \frac{\lambda^3}{2} \int_{\S^2_{+}} \tilde{g}^{0,\lambda}(\vec{W}[f^{0,\lambda},\tilde{g}^{0,\lambda}], X_i)\,
d\mu_{f^{0,\lambda}}\\
& = & \frac{\lambda^3}{2} \int_{\S^2_{+}} W[f^{0,\lambda},\tilde{g}^{0,\lambda}] \xi_i\,d\mu_{f^{0,\lambda}}.
\end{eqnarray*}
The first variation formula for the area yields 
$$
\frac{\partial}{\partial a_i} {\cal A}(\phi^{a,\lambda})|_{a=0} 
= \int_{\S^2_{+}} H[\phi^{0,\lambda}] \varphi_i\,d\mu_{\phi^{0,\lambda}}
+ \int_{\partial \S^2_{+}} g(\tau_i,\eta)\,ds_{\phi^{0,\lambda}}.
$$
Since $g(\tau_i,\eta) \equiv 0$ along $\partial \S^2_{+}$, we get by transforming the integral
$$
\frac{\partial}{\partial a_i} {\cal A}(\phi^{a,\lambda})|_{a=0} = 
\lambda \int_{\S^2_{+}} H[f^{0,\lambda},\tilde{g}^{0,\lambda}] \xi_i\,d\mu_{f^{0,\lambda}}.
$$
But ${\cal A}(\phi^{a,\lambda}) \equiv 2\pi \lambda^2$ for all $a \in S$, therefore we have 
$$
\int_{\S^2_{+}} H[f^{0,\lambda},\tilde{g}^{0,\lambda}] \xi_i\,d\mu_{f^{0,\lambda}} = 0 \quad 
\mbox{ for } i=1,2.
$$
Now if $0 \in S$ is a critical point for $\Wbar(\cdot,\lambda)$, then we also get
$$
\int_{\S^2_{+}} W[f^{0,\lambda},\tilde{g}^{0,\lambda}] \xi_i\,d\mu_{f^{0,\lambda}} = 0
\quad \mbox{ for } i=1,2.
$$
By construction there exist $\alpha,\,\beta_{1,2} \in \R$ such that for 
$\psi_i = \psi_i[w^{0,\lambda},\tilde{g}^{0,\lambda}]$ as in (\ref{eqpsifunctions})
$$
W[f^{0,\lambda},\tilde{g}^{0,\lambda}] = \alpha \psi_0 + \beta_i \psi_i. 
$$
With respect to the metric $g^{0,\lambda} = (f^{0,\lambda})^\ast \tilde{g}^{0,\lambda}$, 
the functions $\xi_i$ are $L^2$-orthogonal to both $W[f^{0,\lambda},\tilde{g}^{0,\lambda}]$
and $H[f^{0,\lambda},\tilde{g}^{0,\lambda}]$. This yields
$$
0 = \big\langle \xi_i, W[f^{0,\lambda},\tilde{g}^{0,\lambda}] \big\rangle_{L^2(\S^2_{+},g^{0,\lambda})}
= \sum_{j=1}^2 \langle \xi_i,\psi_j \rangle_{L^2(\S^2_{+},g^{0,\lambda})} \beta_j
\quad \mbox{ for } i=1,2.
$$
From (\ref{eqtranslation2}), (\ref{eqtranslation3}) we have 
$\|\xi_i - \langle \omega,e_i \rangle\|_{C^0(\S^2_{+})} \leq C \lambda$. 
Recalling that $\psi_i = \sqrt{\frac{3}{2\pi}} \langle \omega,e_i \rangle$ for 
$w = 0$, $\tilde{g} = 0$, we conclude
$$
\Big|\langle \xi_i,\psi_j \rangle_{L^2(\S^2_{+},g^{0,\lambda})} - \sqrt{\frac{2\pi}{3}} \delta_{ij}\Big| \leq C \lambda.
$$
This implies $\beta_1 = \beta_2 = 0$ for $\lambda \leq \lambda_0 = \lambda_0(\Omega)$, and 
we conclude $W[\phi^{0,\lambda}] = \alpha H[\phi^{0,\lambda}]$ as claimed.\\
\\
For the reverse implication assume that $\phi^{0,\lambda}$ is critical for the 
Willmore functional in ${\cal M}^{4,\alpha}_\lambda(S)$, i.e. $W[\phi^{0,\lambda}] = \alpha H[\phi^{0,\lambda}]$
for some $\alpha \in \R$. Then we compute 
\begin{eqnarray*}
\frac{\partial}{\partial a^i}{\cal W}(\phi^{a,\lambda})|_{a= 0} & = &
\langle \vec{W}[\phi^{0,\lambda}], \frac{\partial \phi^{a,\lambda}}{\partial a^i}|_{a=0} \rangle_{L^2}\\
& = & \alpha \langle \vec{H}[\phi^{0,\lambda}], \frac{\partial \phi^{a,\lambda}}{\partial a^i}|_{a=0} \rangle_{L^2}\\
& = & - \frac{\partial}{\partial a^i}{\cal A}(\phi^{a,\lambda})|_{a= 0}\\
& = & 0.
\end{eqnarray*}
Hence $a=0$ is a critical point of $\Wbar(\cdot,\lambda)$, which finishes the proof 
of the theorem. \kasten

We finally prove a purely local existence result. 

\begin{theorem} \label{thmlocal} Let $\Omega$ be a bounded domain of class $C^{12}$.
If $a \in S = \partial \Omega$ is a nondegenerate critical point of $H^S$,
then there exists a $C^1$ curve $\gamma(\lambda) \in S$ for 
$\lambda \in [0,\lambda_0)$, such that $\gamma(0) = a$ and each 
$\phi^{\gamma(\lambda),\lambda}$, $\lambda > 0$, is a critical point of 
${\cal W}(f)$ in ${\cal M}^{4,\alpha}_\lambda(S)$.
\end{theorem}

We need the following calculus lemma.

\begin{lemma} Let $u \in C^2\big(S \times (-\lambda_0,\lambda_0)\big)$ be a given function
satisfying $u(\cdot,0) \equiv 0$, and let $v: S \times (-\lambda_0,\lambda_0) \to \R$ be
defined by
$$
v(a,\lambda) = \begin{cases} 
\lambda^{-1} u(a,\lambda) & \mbox{ for }\lambda \neq 0,\\
\partial_\lambda u(a,0) & \mbox{ for } \lambda = 0,
\end{cases}
$$
Then $v$ is of class $C^1\big(S \times (-\lambda_0,\lambda_0)\big)$, having the derivatives
\begin{eqnarray*}
\nabla v (a,\lambda) & = & \begin{cases}
\lambda^{-1} \nabla u(a,\lambda) & \mbox{ for }\lambda \neq 0,\\
\partial_\lambda \nabla u(a,0) & \mbox{ for } \lambda = 0,
\end{cases}\\
\partial_\lambda v(a,\lambda) & = & 
\begin{cases}
\lambda^{-2} \big(\lambda \partial_\lambda u(a,\lambda) - u(a,\lambda)\big) & \mbox{ for } \lambda \neq 0,\\
\frac{1}{2}\, \partial_\lambda^2 u(a,0) & \mbox{ for } \lambda = 0.
\end{cases}
\end{eqnarray*}
\end{lemma}

{\em Proof. } We have using $u(x,0) = 0$
$$
|\lambda^{-1} u(x,\lambda) - \partial_\lambda u(a,0)| 
= \big|\mint_0^\lambda (\partial_\lambda u(x,s) - \partial_\lambda u(a,0))\,ds\big|
\to 0 \quad \mbox{ for } x \to a,\,\lambda \to 0.
$$
This shows that $v$ is continuous. For the $C^1$ property it is sufficient to
prove that the stated derivatives are also continuous. In the case of $\nabla v$
the argument above applies (noting that $\nabla u$ is $C^1$ by assumption).
For $\partial_\lambda v$ we compute
\begin{eqnarray*}
&&\lambda^{-2} \big(\lambda \partial_\lambda u(x,\lambda) - u(x,\lambda)\big) - \frac{1}{2} \partial_\lambda^2 u(a,0)\\ 
& = & \lambda^{-2} \int_0^\lambda \big(\partial_\lambda u(x,\lambda) - \partial_\lambda u(x,t)\big)\,dt 
- \frac{1}{2} \partial_\lambda^2 u(a,0)\\
& = & \lambda^{-2} \int_0^\lambda \int_t^\lambda \big(\partial_\lambda^2 u(x,s) - \partial_\lambda^2 u(a,0)\big)\,ds dt\\
& \to & 0 \quad \mbox{ for } a \to x,\,\lambda \to 0.
\end{eqnarray*}
The lemma is proved. \kasten

{\em Proof of Theorem. } We apply the lemma to the function $u(a,\lambda) = \nabla \bar{{\cal W}}(a,\lambda)$,
where $\bar{{\cal W}}(\cdot,0) \equiv 2\pi$ and hence $u(a,0) = \nabla \bar{{\cal W}}(a,0) \equiv 0$. 
This needs $\bar{{\cal W}} \in C^3\big(S \times (-\lambda_0,\lambda_0)\big)$, which is true 
for $\Omega \in C^{12}$. We obtain from (\ref{eqmixedderivatives}) 
(taking one more derivative $\nabla_a$)
\begin{eqnarray*}
v(a,0) & = & \partial_\lambda \nabla \bar{{\cal W}}(a,0) = - \pi \nabla H^S(a)\\
\nabla v(a,0) & = & \partial_\lambda \nabla^2 \bar{{\cal W}}(a,0) = - \pi \nabla^2 H^S(a).
\end{eqnarray*}
Now assume for $0 \in S$ that $\nabla H^S(0) = 0$ and $\nabla^2 H^S(0)$ nondegenerate.
Then the implicit function theorem, applied to $v(a,\lambda)$, yields a neigborhood
$U \times (-\ve,\ve)$ and a $C^1$-curve $a = \gamma(\lambda)$, such that for
$(a,\lambda) \in U \times (-\ve,\ve)$ one has
$$
v(a,\lambda) = 0 \quad \Leftrightarrow \quad a = \gamma(\lambda).
$$
For $\lambda \neq 0$ we thus get
$$
\nabla \bar{{\cal W}}(a,\lambda) = 0 \quad \Leftrightarrow \quad a = \gamma(\lambda).
$$
The theorem now follows from Theorem \ref{thmreduction}. \kasten

\section{Appendix: Construction of the barycenter}
The concept of Riemannian barycenter is due to Karcher \cite{Kar77}. For 
our purposes we only need a local version, which does not involve e.g. 
Riemannian comparison theory. 
Let $U = D_\delta(0) \subset \R^2$, $V = B_{\frac{3}{2}}(0) \subset \R^3$. 
For $x \in U$, $v \in V$ we put
$$
c_{x,v}:[0,1] \to Z_2,\,c_{x,v}(t) = x + t v.
$$
Further let $X = \{\phi \in C^2([0,1],\R^3): \phi(0) = \phi'(0) = 0\}$ and 
$$
X_\ve = \{\phi \in X: \|\phi\|_{C^0([0,1])} < \ve\}.
$$ 
We finally put $G_\ve = \{\tilde{g} \in C^l(\bar{Z_2},\R^{3 \times 3}):
\|\tilde{g} - \delta\|_{C^l(Z_2)} < \ve\}$ for $l \geq 1$, and consider
$$
F:U \times V \times X_\ve \times G_\ve \to C^0([0,1],\R^3),\,
F[x,v,\phi,\tilde{g}] = c'' + \tilde{\Gamma} \circ c (c',c')|_{c = c_{x,v} +\phi}.
$$
We claim that $F$ is of class $C^{l-1}$. Write $F = F_2 \circ F_1$
where $F_1$ is the affine map
$$
F_1:U \times V \times X_\ve \to C^2([0,1],\R^3),\,F_1[x,v,\phi] = c_{x,v} + \phi.
$$
$F_1$ is continuous and hence smooth. 
The nonlinear map $F_2$ is given by 
$$
F_2:C^2([0,1],Z_2) \times G_\ve \to C^0([0,1],\R^3),\,
F_2[c,\tilde{g}] = c'' + \tilde{\Gamma} \circ c (c',c').
$$
The composition $C^2 \times C^{l-1} \to C^0$, $(c,\tilde{\Gamma}) \mapsto \tilde{\Gamma} \circ c$,
is of class $C^{l-1}$. Namely differentiating $l-1$ times with respect to $c$ leaves 
exactly a $C^0$ function. Since we can build $F_2$ from $\tilde{\Gamma}\circ c$
by linear or bilinear operations, it is also of class $C^{l-1}$. Assuming 
from now on $l \geq 2$, we have 
\begin{eqnarray*}
D_c F_2[c,\tilde{g}] \phi & = & \phi'' 
+ 2\, \tilde{\Gamma} \circ c(\phi',c') + (D\tilde{\Gamma})\circ c (\phi,c',c')\\
D_{\tilde{g}} F_2[c,\tilde{g}] h & = & 
\frac{1}{2} \Big(\tilde{g}^{kp}(\partial_i h_{jp} + \partial_j h_{ip} - \partial_p h_{ij})\Big) \circ c\, (c^i)' (c^j)' e_k\\ 
&& - \frac{1}{2} \Big(\tilde{g}^{km} h_{mn} \tilde{g}^{np} 
\left(\partial_i \tilde{g}_{jp} + \partial_j \tilde{g}_{ip} - \partial_p \tilde{g}_{ij}\right)\Big) \circ c\, (c^i)' (c^j)' e_k.
\end{eqnarray*}
In particular 
$$
F[x,v,0,\delta] = 0 \quad \mbox{ and } \quad
D_\phi F[x,v,0,\delta] \psi = \psi''.
$$
The map $D_\phi F[x,v,0,\delta]:X \to C^0([0,1],\R^3)$ is an isomorphism, in fact 
the equation $\psi'' = f$ has the unique solution $\psi \in X$ given by 
$$
\psi(u) = \int_0^u \int_0^s f(t)dt\,ds. 
$$
By the implicit function theorem, the set of solutions of $F[x,v,\phi,\tilde{g}] = 0$ 
near $[0,v_0,0,\delta]$ is given as a $C^{l-1}$ graph
$$
{\bf \phi} = {\bf \phi}[x,v,\tilde{g}],
$$
i.e. the corresponding curves ${\bf c}[x,v,\tilde{g}] = c_{x,v} + {\phi}[x,v,\tilde{g}]$ are 
geodesics with respect to $\tilde{g}$ having initial data $c(0) = x$, $c'(0) = v$. The 
exponential mapping is now given by 
$$
\exp:U \times V \times G_\ve \to Z_2,\,\exp_x^{\tilde{g}}(v) = {\bf c}[x,v,\tilde{g}](1).
$$
Now $D\exp_x^{\delta} = {\rm Id}_{\R^3}$, Thus for $l \geq 2$ and $\ve > 0$ small we get 
$$
\|D\exp^{\tilde{g}}_x - {\rm Id}_{\R^3}\|_{C^0(V)} \leq \ve_0 \quad 
\mbox{ for } \tilde{g} \in G_\ve,\,x \in U = D_\delta(0).
$$
This gives for $v,w \in V$
\begin{eqnarray*}
|\exp^{\tilde{g}}_x (v) -  \exp^{\tilde{g}}_x (w)| & = & 
\Big|\int_0^1 D\exp^{\tilde{g}}_x((1-t)w + tv)\cdot(v-w)\,dt\Big|\\
& \geq & |v-w| - \Big|\int_0^1 \big(D\exp^{\tilde{g}}_x((1-t)w + tv) - {\rm Id}_{\R^3}\big)\cdot(v-w)\,dt\Big|\\
& \geq & (1-\ve_0) |v-w|. 
\end{eqnarray*}
This shows that $\exp_x^{\tilde{g}}$ is injective on $V = B_{\frac{3}{2}}(0)$. 
We further estimate 
$$
|\exp^{\tilde{g}}_x(v) - v| =  
\Big|\int_0^1 \big(D\exp^{\tilde{g}}_x(tv) - {\rm Id}_{\R^3}\big)\,dt \cdot v\Big| 
\leq \ve_0\,|v|.
$$
We now show that $\exp_x^{\tilde{g}}(V) \cap B_{\frac{5}{4}}(0)$ is a closed
subset of $B_{\frac{5}{4}}(0)$. Assume that $\exp_x^{\tilde{g}}(v_k) \to p \in B_{\frac{5}{4}}(0)$.
From the above we then have 
$$
|v_k| - \frac{5}{4} < |v_k| - |\exp_x^{\tilde{g}}(v_k)| \leq |v_k - \exp_x^{\tilde{g}}(v_k)| \leq \ve_0 |v_k|,
$$
which implies $|v_k| \leq (1-\ve_0)^{-1}\, \frac{5}{4} < \frac{3}{2}$ for appropriate $\ve_ 0 > 0$.
Up to a subsequence, we thus have $v_k \to v \in V$ and $\exp_x^{\tilde{g}}(v) = p$. Now 
$\exp_x^{\tilde{g}}(V) \cap B_{\frac{5}{4}}(0)$ is also open by the inverse function theorem, 
hence we have $B_{\frac{5}{4}}(0) \subset \exp_x^{\tilde{g}}(V)$, and we obtain the inverse
$$
(\exp_x^{\tilde{g}})^{-1}: B_{\frac{5}{4}}(0) \to V.
$$
Of course we are not claiming that $\exp_x^{\tilde{g}}$ maps all of $V$ into $B_{\frac{5}{4}}(0)$. 
The inverse is of class $C^{l-1}$ in all variables 
$x \in U$, $p \in  B_{\frac{5}{4}}(0)$ and $\tilde{g} \in C^l(Z_2)$. Namely let 
$\exp_{x_0}^{\tilde{g}_0}(v_0) = p_0 \in B_{\frac{5}{4}}(0)$, where $v_0 \in V$.
Consider the equation
$$
\exp_x^{\tilde{g}}(v) - p = 0.
$$
By the implicit function theorem, the set of solutions has a local 
representation $v = v[x,p,\tilde{g}]$ which is of class $C^{l-1}$.
But the local inverse equals the global inverse, and hence also 
the global inverse is of class $C^{l-1}$ as claimed.

\begin{lemma}[two-dimensional barycenter] \label{lemmabarycenter} 
Assume $w:\S^2_{+} \to \R$, $\tilde{g}:Z_2 \to \R^{3 \times 3}$ 
belong to the neighborhoods $W_\ve$, $G_\ve$ given by 
$$
\|w\|_{C^1(\S^2_{+})} < \ve \quad \mbox{ and } \quad 
\|\tilde{g}-\delta\|_{C^l(Z_2)} < \ve \quad \mbox{ where } l\geq 2.
$$
For $\ve > 0$ small we then have a welldefined function
$$
X[w,\tilde{g}]:U \to \R^2,\,X[w,\tilde{g}](x) 
= -\pi_{\R^2}  \Big(\int_{\S^2_{+}} (\exp^{\tilde{g}}_x)^{-1} (f(\omega))\,d\mu_g(\omega)\Big),
$$
and there is a unique point $x \in U$ with $X[w,\tilde{g}](x) = 0$. 
This point $x =  C[w,\tilde{g}]$ is called the two-dimensional barycenter
of (the radial graph of) $w$ with respect to $\tilde{g}$.
The map $C[w,\tilde{g}]$ is of class $C^{l-1}$.
\end{lemma}

{\em Proof. }Let $f(\omega) = \omega + w(\omega)$. Fixing a coordinate system 
on $\S^2_{+}$, we consider the map 
\begin{equation}
\label{bary1}
U \times W_\ve \times G_\ve \to C^0(\S^2_{+},\R^3),\,[x,w,\tilde{g}] 
\mapsto (\exp^{\tilde{g}}_x)^{-1} \circ f\, \sqrt{\det g}. 
\end{equation}
By standard rules for product and composition, the right hand side 
belongs to $C^0(\S^2_{+},\R^3)$; in particular $X[w,\tilde{g}]$ is well-defined.
We claim that the map (\ref{bary1}) is of class $C^{l-1}$ in all three variables. 
For this we recall that $\Psi[x,p,\tilde{g}] = (\exp_x^{\tilde{g}})^{-1}(p)$
is of class $C^{l-1}$. For $\omega \in \S^2_{+}$ fixed we have the $C^{l-1}$ composition
$$
\begin{array}{ccccc}
U \times W_\ve \times G_\ve & \stackrel{C^\infty}{\rightarrow} & U \times B_{\frac{5}{4}}(0) \times G_\ve 
& \stackrel{\Psi}{\rightarrow} & V\\
(x,w,\tilde{g}) & \mapsto & (x,f(\omega),\tilde{g}) 
&\mapsto & \Psi[x,f(\omega),\tilde{g}].
\end{array}
$$
Now all derivatives with respect to $x,w,\tilde{g}$ up to order $l-1$
depend also continuously on $\omega$, which yields the claim.
For $\tilde{g} = \delta$ we have $(\exp_x)^{-1}(p) = p-x$ which implies
$$
X[w,\delta](x) = \mu_g(\S^2_{+})\,
\Big(x - \pi_{\R^2} \mint_{\S^2_{+}} f(\omega)\,d\mu_g(\omega)\Big),
$$
in particular $X[0,\delta](0) = 0$ and $D_x X[0,\delta](x) = 2\pi\, {\rm Id}_{\R^2}$.
Thus by the implicit function theorem there is a unique point $x \in U$ with 
$X[w,\tilde{g}](x) = 0$, and the resulting map $x = C[w,\tilde{g}]$ is of class
$C^{l-1}$. \kasten

From the proof we note the explicit formula 
\begin{equation}
\label{eqbary2}
C[w,\delta] = \pi_{\R^2} \Big(\mint_{\S^2_{+}} f(\omega)\,d\mu_g(\omega)\Big).
\end{equation}
We consider the two coordinates $C^i[f,\tilde{g}]$ of the barycenter as functionals 
depending on $w$ resp. $f$, and we now compute the corresponding $L^2$ gradient. 
Consider a compactly supported variation of $f$ in direction $\phi = \varphi \nu$. 
Then we have 
$$
\frac{\partial}{\partial \ve} (g_{\ve})_{ij}|_{\ve = 0} = - 2\varphi h_{ij} 
\quad \mbox{ and } \quad
\frac{\partial}{\partial \ve} d\mu_{g_\ve}|_{\ve = 0} = - \varphi H \,d\mu_g. 
$$
The first variation of $X[f,\tilde{g}]$ is then 
\begin{eqnarray*}
\frac{\partial}{\partial \ve} X[f_\ve,\tilde{g}](x)|_{\ve =0} & = & 
- \pi_{\R^2} \int_{\S^2_{+}} D ((\exp^{\tilde{g}}_x)^{-1})(f(\omega)) \cdot \phi(\omega)\,d\mu_g(\omega)\\
&& + \pi_{\R^2} \int_{\S^2_{+}} (\exp^{\tilde{g}}_x)^{-1}(f(\omega)) H(\omega) \varphi(\omega)\,d\mu_g(\omega).
\end{eqnarray*}
By definition of the barycenter we have 
\begin{eqnarray*}
0 & = & \frac{\partial}{\partial \ve} X\big[f_\ve,\tilde{g}](C[f_\ve,\tilde{g}])|_{\ve =0}\\
& = &  -\pi_{\R^2} \int_{\S^2_{+}} 
D\big((\exp^{\tilde{g}}_x)^{-1}\big)(f(\omega)) \phi(\omega)\,d\mu_g(\omega)|_{x = C[f,\tilde{g}]}\\
&& + \pi_{\R^2} \int_{\S^2_{+}} 
(\exp^{\tilde{g}}_x)^{-1}(f(\omega)) H(\omega) \varphi(\omega)\,d\mu_g(\omega)|_{x = C[f,\tilde{g}]}\\
&& - \pi_{\R^2} \int_{\S^2_{+}} 
D_x (\exp^{\tilde{g}}_x)^{-1} (f(\omega))\,d\mu_g(\omega)|_{x = C[f,\tilde{g}]} 
\cdot \frac{\partial}{\partial \ve} C[f_\ve,\tilde{g}]|_{\ve =0}.
\end{eqnarray*}
This implies the formula 
\begin{eqnarray*}
\frac{\partial}{\partial \ve}C[f_\ve,\tilde{g}]|_{\ve =0} & = & \Big(\pi_{\R^2} \int_{\S^2_{+}}
D_x (\exp^{\tilde{g}}_x)^{-1}(f(\omega))\,d\mu_g(\omega)\Big)^{-1}|_{x = C[f,\tilde{g}]} \\
&\cdot& \Big(-\pi_{\R^2} \int_{\S^2_{+}} 
D\big((\exp^{\tilde{g}}_x)^{-1}\big)(f(\omega)) \phi(\omega)\,d\mu_g(\omega)|_{x = C[f,\tilde{g}]}\\
&& + \pi_{\R^2} \int_{\S^2_{+}} 
(\exp^{\tilde{g}}_x)^{-1}(f(\omega)) H(\omega) \varphi(\omega)\,d\mu_g(\omega)|_{x = C[f,\tilde{g}]}\Big).
\end{eqnarray*}
Under reparametrizations of $f$ the barycenter remains the same, hence 
the $L^2$ gradient of $C^i[f,\tilde{g}]$ is normal along $f$. Taking 
the $\tilde{g}$ inner product with $\nu$ yields a scalar function, 
which we denote by ${\rm grad}_{L^2}\,C^i[w,\tilde{g}]$ in slight 
abuse of notation. We now conclude
\begin{align} \label{eqbary3}
\sum_{i=1}^2 {\rm grad}_{L^2}\,C^i[f,\tilde{g}] e_i  = &  \nonumber \Big(\pi_{\R^2} \int_{\S^2_{+}}
D_x (\exp_x^{\tilde{g}})^{-1} (f(\omega))\,d\mu_g(\omega)\Big)^{-1}|_{x = C[f,\tilde{g}]}\\ 
& \cdot \pi_{\R^2} \Big(
(\exp_x^{\tilde{g}})^{-1}(f) H - D\big((\exp_x^{\tilde{g}})^{-1}\big)(f) \nu \Big)|_{x = C[f,\tilde{g}]}.
\end{align}
In the Euclidean case $\tilde{g} = \delta$ we have $\exp_x v = x+v$, which yields for $i = 1,2$
\begin{eqnarray*}
\pi_{\R^2} \int_{\S^2_{+}} D_x (\exp_x)^{-1}(f(\omega))\,d\mu_g(\omega) & = & - \mu_g(\S^2_{+})\, {\rm Id}_{\R^2},\\
{\rm grad}_{L^2}\,C^i[f,\delta] & = & \frac{1}{\mu_g(\S^2_{+})} \big\langle \nu - (f - C[f,\delta]) H, e_i \big\rangle.
\end{eqnarray*}
Specializing further to $f_0(\omega) = \omega$, we see 
\begin{equation}
\label{eqbary4}
{\rm grad}_{L^2}\,C^i[f_0,\delta](\omega) = - \frac{3}{2\pi} \langle \omega,e_i \rangle.
\end{equation}
For $w \in C^{k,\alpha}(\S^2_{+})$ and $\tilde{g} \in C^l(\overline{Z}_2,\R^{3 \times 3})$ 
where $l \geq k+1$, one deduces ${\rm grad}_{L^2} C^i[w,\tilde{g}] \in C^{k-2,\alpha}(\S^2_{+})$.
Moreover as a functional into $C^{k-4,\alpha}(\S^2_{+})$, it is of class $C^{l-k+1}$.

\begin{minipage}[t]{8cm}
Roberta Alessandroni\\
Mathematisches Institut\\
Albert-Ludwigs-Universit\"at Freiburg\\
Eckerstra{\ss}e 1, D 79104 Freiburg\\
roberta.alessandroni@math.uni-freiburg.de
\end{minipage}
\ \
\begin{minipage}[t]{7cm}
Ernst Kuwert\\
Mathematisches Institut\\
Albert-Ludwigs-Universit\"at Freiburg\\
Eckerstra{\ss}e 1, D 79104 Freiburg\\
ernst.kuwert@math.uni-freiburg.de
\end{minipage}

\end{document}